\documentclass{amsart}

\usepackage{graphicx}
\usepackage{natbib}
\usepackage{amsmath}

\newtheorem{remark}{Remark}
\newtheorem{notation}{Notation}

\title[Long-waves models in presence of vorticity]{Fully nonlinear long-waves models in presence of vorticity}

\author[A. Castro and  D. Lannes]%
{Angel Castro,
David Lannes
}
\thanks{angel.castro@uam.es, David.Lannes@ens.fr}
\address{$^1$ Departamento de Matem\'aticas UAM, Instituto de Ciencias Matem\'aticas CSIC, Campus de Cantoblanco, 28049  Madrid, Spain\\
$^2$DMA, Ecole Normale Sup\'erieure et CNRS UMR 8553, 45 rue
  d'Ulm, 75005 Paris, France}

\newcommand{\abs}[1]{\vert#1\vert}
\newcommand{\babs}[1]{\big\vert#1\big\vert}

\newcommand{\dsp}{\displaystyle}
\newcommand{\bU}{{\mathbf U}}

\newcommand{\ovV}{\overline{V}}
\newcommand{\ovv}{\overline{v}}

\newcommand{\uU}{\underline{U}}

\newcommand{\bom}{{\boldsymbol \omega}}

\newcommand{\uV}{{\underline V}}

\newcommand{\uw}{{\underline w}}

\newcommand{\R}{{\mathbb R}}

\newcommand{\dt}{\partial_t}
\newcommand{\dz}{\partial_z}

\newcommand{\curlm}{\mbox{\textnormal{curl}}^\mu}

\newcommand{\dive}{\mbox{\textnormal{div} }}

\newcommand{\divem}{\mbox{\textnormal{div}}^\mu}
\newcommand{\surf}{{\vert_{\rm surf}}}
\newcommand{\bott}{{\vert_{\rm bott}}}

\newcommand{\eps}{\varepsilon}

\newcommand{\dx}{\partial_x}

\newcommand{\ep}{\varepsilon}

\newcommand{\pa}{\partial}

\newcommand{\vs}{V_{\rm sh}^*}
\begin{document}

\maketitle

\begin{abstract}
We study here Green-Naghdi type equations (also called fully nonlinear
Boussinesq, or Serre equations) modeling the propagation of large
amplitude waves in shallow water. The novelty here is that we allow
for a general vorticity, hereby allowing complex interactions between
surface waves and currents. We show that the {\it a priori}
$2+1$-dimensional dynamics of the vorticity can be reduced to a
finite cascade of two-dimensional equations: with a mechanism
reminiscent of turbulence theory, vorticity effects contribute to the
averaged momentum equation through a Reynolds-like tensor that can be
determined by a cascade of equations. Closure is obtained at the
precision of the model at the second order of this cascade. We also
show how to reconstruct the velocity field in the $2+1$ dimensional
fluid domain from this set of $2$-dimensional equations and exhibit
transfer mechanisms between the horizontal and vertical components of
the vorticity, thus opening perspectives for the study of rip currents
for instance.
\end{abstract}

\noindent
{\bf Keywords.}
Green-Naghdi equation, Serre equations, Fully nonlinear
  Boussinesq equations, Shallow water, Rotational flows,
Wave-current interactions, Rip-currents

\section{Introduction}

The equations describing the motion of an inviscid and
incompressible fluid of constant density $\rho$ and
delimited from a above by a free surface $\{z=\zeta(t,X)\}$ ($X\in
\R^2$) and below by a non moving bottom $\{z=-H_0+b(X)\}$ are given by
the so called free surface Euler equations. Denoting by $\bU=({\bf
  V}^T,{\bf w})^T$ and $P$
the velocity and pressure fields, these
equations can be written
\begin{eqnarray}\label{Eulereq}
\dt \bU+\bU\cdot \nabla_{X,z}\bU&=&-\frac{1}{\rho}\nabla_{X,z}P -g {\bf e}_z,\\
\label{incomp}
\nabla_{X,z}\cdot \bU&=&0,
\end{eqnarray}
in the fluid domain $\Omega_t=\{(X,z)\in \R^{2+1},
-H_0+b(X)<z<\zeta(t,X)\}$; they are complemented with the boundary conditions
\begin{eqnarray}
\label{kineticeq}
\dt \zeta -\bU_{\surf}\cdot N&=&0 \quad(\mbox{with }N=(-\nabla\zeta^T,1)^T),\\
P_\surf&=&\mbox{constant}
\end{eqnarray}
at the surface, and
\begin{equation}
\label{condbott}
\bU_\bott\cdot N_b=0\quad (\mbox{with }N_b=(-\nabla b^T,1)^T))
\end{equation}
at the bottom.

It is well-known that the kinetic equation (\ref{kineticeq}) can be
restated as a mass conservation equation
\begin{equation}\label{masscons}
\dt \zeta+\nabla\cdot (h\ovV)=0,
\end{equation}
where $h$ is the total depth of the fluid and $\ovV$ the vertical
average of the horizontal component ${\bf V}$ of the velocity (equivalently, $h\ovV$
is the total discharge),
\begin{equation}\label{defhovV}
h(t,X):=H_0+\zeta(t,X)-b(X),\qquad
\ovV(t,X)=\frac{1}{h}\int_{-H_0+b(X)}^{\zeta(t,X)} {\bf V}(t,X,z)dz.
\end{equation}
\begin{notation}
We decompose any function $f$ defined on
  $\Omega_t$ as an averaged part and a zero mean component, using the notation
$$
f(t,X,z)=\overline{f}(t,X)+f^*(t,X,z),\quad\mbox{with
}\quad \overline{f}=\frac{1}{h}\int_{-H_0+b(X)}^{\zeta(t,X)}f(t,X,z)dz
$$
and $f^*=f-\overline{f}$.
\end{notation}
It is therefore quite natural to look for another equation that would
complement (\ref{masscons}) to form a closed system of two evolution
equations on $\zeta$ and $\ovV$. Decomposing the horizontal velocity
field as
\begin{equation}
\label{decompV}
{\bf V}(t,X,z)=\ovV(t,X)+V^*(t,X,z),
\end{equation}
and integrating vertically the horizontal component of
(\ref{Eulereq}), one obtains classically
$$
\dt (h\ovV)+\nabla\cdot ( \int_{-H_0+b}^\zeta {\bf V}\otimes {\bf
  V})+\int_{-H_0+b}^\zeta \nabla P=0.
$$

Since, by construction, the vertical average of $V^*$ vanishes, we
finally obtain as in \cite{Teshukov} the following set of evolution equations on $\zeta$ and
$\ovV$,
\begin{equation}\label{Eulerav}
\begin{cases}
\dsp \dt \zeta+\nabla\cdot (h\ovV)=0,\\
\dsp \dt (h\ovV)+\nabla\cdot (h\ovV\otimes \ovV)+\nabla\cdot ( \int_{-H_0+b}^\zeta V^*\otimes V^*)+\int_{-H_0+b}^\zeta \nabla P=0.
\end{cases}
\end{equation}
We shall refer to (\ref{Eulerav}) as the {\it averaged Euler
  equations}. These equations are exact but too complex to work with
(because $V^*$ and $\nabla P$ are not closed expressions of $\zeta$
and $\ovV$) and they are replaced by simpler approximate equations for
practical purposes; we shall consider here approximations of these equations in {\it shallow
water}, i.e.  when the depth is small compared to the typical
horizontal length.

Let us first consider the case of {\it irrotational flows} for which
(\ref{Eulereq}) and (\ref{incomp}) are complemented by the condition
\begin{equation}
\label{irrot}
\nabla_{X,z}\times\bU=0.
\end{equation}
 In the shallow water regime, that is, when $\mu:=H_0^2/L^2\ll 1$ (with
$L$ the typical horizontal scale), it is well-known that the flow is
columnar at leading order in the sense that the horizontal velocity
does not depend at leading order on the vertical variable $z$. The
{\it ``Reynolds'' tensor}
$$
{\bf R}:=\int_{-1+ b}^{ \zeta} V^*\otimes V^*
$$
is therefore a second order term (it is of size $O(\mu^2)$).
The terminology ``Reynolds tensor''
is {\it stricto sensu} improper here since space derivatives do not
commute with averaging (consisting in vertical integration here); we
however use it because the analogy with Reynolds turbulence theory is instructive.\\
 It is
also classical in shallow water that the pressure is hydrostatic at
leading order, $P(t,X,z)= \rho g (\zeta(t,X)-z)+O(\mu)$. At leading
order, the averaged Euler equations (\ref{Eulerav}) are therefore
formally approximated by the nonlinear shallow water (or Saint-Venant)
equations,
\begin{equation}\label{StVenant}
\begin{cases}
\dsp \dt \zeta+\nabla\cdot (h\ovV)=0,\\
\dsp \dt (h\ovV)+ gh \nabla\zeta+\nabla\cdot (h\ovV\otimes \ovV)=0
\end{cases}
\end{equation}
(see \cite{Ovs,KN,AlvarezLannes,Iguchi} for a justification of this
approximation). This model is widely used but misses for instance
dispersive effects that can be very important in coastal
oceanography. This is the reason why a more precise model taking into
account the first order terms (with respect to $\mu$) is used for
applications. As already said, the Reynolds-like tensor is a second order
term and can still be neglected at this level of approximation, but
non-hydrostatic components of the pressure must be taken into
account. The resulting equations are known as the Green-Naghdi (or
Serre, or fully nonlinear Boussinesq equations, see \cite{LannesBonneton,Betal} for
recent reviews). Under the formulation derived in \cite{BCLMT}, these
equations can be written
\begin{equation}\label{Green-Naghdi}
\begin{cases}
\dsp \dt \zeta+\nabla\cdot (h\ovV)=0,\\
\dsp (I+h{\mathcal T}\frac{1}{h})\big(\dt (h\ovV)+\nabla\cdot
(h\ovV\otimes \ovV)\big)+g h \nabla\zeta+h{\mathcal Q_1}(\ovV)=0,
\end{cases}
\end{equation}
where the non-hydrostatic effects are taken into account through the
operators ${\mathcal T}$ and ${\mathcal Q}(\cdot)$ defined as
\begin{equation}\label{defTintro}
{\mathcal T}V=-\frac{1}{3h}\nabla(h^3\nabla\cdot
V)+\frac{1}{2h}\big[\nabla (h^2\nabla b\cdot V)-h^2\nabla b\nabla\cdot
V\big]+\nabla b \nabla b\cdot V
\end{equation}
and, writing $V^\perp=(-V_2,V_1)^\perp$,
\begin{equation}\label{defQ1intro}
{\mathcal Q}_1(V)=-2{\mathcal R}_1\big(\partial_1V\cdot \partial_2
V^\perp+(\nabla\cdot V)^2\big)+{\mathcal R}_2\big(V\cdot (V\cdot
\nabla)\nabla b\big)
\end{equation}
with
$$
{\mathcal R}_1w=-\frac{1}{3h}\nabla(h^3w)-\frac{h}{2}w\nabla b,\qquad
{\mathcal R}_2w=\frac{1}{2h}\nabla(h^2w)+w\nabla b.
$$
It is possible to derive from (\ref{Green-Naghdi}) an equation for the
local conservation of energy, namely,
\begin{equation}\label{NRJintro}
\dt {\mathfrak e}+\nabla\cdot {\mathfrak F}=0
\end{equation}
where the energy ${\mathfrak e}$ is the sum of the potential energy
${\mathfrak e}_p$ and of the kinetic energy ${\mathfrak e}_k$ given by
$$
{\mathfrak e}_p=\frac{1}{2}g \zeta^2,
\qquad
{\mathfrak e}_k=\frac{1}{2}h \abs{\ovV}^2+\frac{1}{2} h \big(\frac{1}{3}\babs{h\nabla\cdot
  \ovV-\frac{3}{2}\nabla b \cdot \ovV}^2+\frac{1}{4}\abs{\nabla b\cdot
\ovV}^2\big),
$$
and where the flux ${\mathfrak F}$ is given by the expression
$$
{\mathfrak F}=(g\zeta h+{\mathfrak e}_k+{\mathfrak q})\ovV,
$$
with
$$
{\mathfrak q}=-\frac{1}{3}h^2(\dt +\ovV\cdot \nabla)(h\nabla\cdot \ovV)+\frac{1}{2}h^2(\pa_t +\ovV\cdot\nabla)(\nabla b\cdot \ovV).
$$

The Green-Naghdi equations have been rigorously justified in
\cite{Makarenko,Li,AlvarezLannes} (see also the monograph
\cite{Lannes_book} and references therein). The Green-Naghdi system is now
the most popular model for the numerical simulation of coastal flows,
even in configurations that include vanishing depth (shoreline) and
wave breaking (see for instance \cite{CKDKC,Cienfuegos,LGH,BCLMT,Delis,Dutykh,Mario,MarcheLannes}).

Despite their many advantages, the Green-Naghdi equations
(\ref{Green-Naghdi}) can only account for configurations where
rotational effects are absent (i.e. when the assumption (\ref{irrot})
holds). This is unfortunately not the case when waves propagate in a
zone where currents are present, or when vorticity is created by
anisotropic dissipation due to wave breaking as for rip-currents (see for
instance \cite{Hammack,Chen1}). More generally, the full coupling between
currents and surface waves is still largely not understood. The presence of a non zero vorticity makes
the analysis more difficult. Indeed, $d+1$ dimensional irrotational flows ($d$
being the horizontal dimension) are $d$ dimensional in
nature: as shown by Zakharov (\cite{Zakharov}), the full Euler equations
can be reduced to an Hamiltonian formulation in terms of $\zeta(t,X)$ and
$\psi(t,X)=\Phi(t,X,\zeta(t,X))$ where $\Phi$ is a scalar velocity
potential (i.e. $\nabla_{X,z}\Phi=\bU$). Both $\zeta$ and $\psi$ are
independent of the vertical variable and Zakharov's Hamiltonian
formulation of the full water waves problem is therefore not qualitatively different in this aspect
from averaged models such as the Saint-Venant or Green-Naghdi equations
(\ref{StVenant}) and (\ref{Green-Naghdi}). In the rotational setting,
the picture is drastically different, and $d+1$ dimensional flows are
truly $d+1$ dimensional. A Hamiltonian formulation of the water waves
equations in presence of vorticity generalizing Zakharov's formulation
has recently been derived in \cite{CastroLannes}. The evolution of
$\zeta$ and $\psi$ must then be coupled to the evolution equation on
the vorticity $\bom=\nabla_{X,z}\times \bU$,
\begin{equation}\label{intreqvort}
\dt \bom+\bU\cdot \nabla_{X,z}\bU=\bU\cdot \nabla_{X,z}\bom \quad
\mbox{ in }\quad \Omega_t,
\end{equation}
which is $d+1$ dimensional. The reduction to a $d$ dimensional set of
equations as the Saint-Venant or Green-Naghdi equations is therefore a
qualitative jump and is not  {\it a priori}
obvious. Technically, in the rotational setting, the ``Reynolds'' tensor
${\mathbf R}$ is no longer a second order term and contributes to the
momentum equation through a coupling with the $d+1$ dimensional
dynamics of the vorticity equation.

 Several
approaches have been proposed to get around this difficulty. Following
 \cite{Bowen} and  \cite{LH}, many models use
a time-averaging approach where radiation stresses due to the short
wave motion are considered as a forcing term in the momentum equations
(see \cite{SvendsenPutrevu} for a review). A generalization of the
Green-Naghdi/fully nonlinear Boussinesq/Serre equations
(\ref{Green-Naghdi}) able to handle the presence of vorticity would be
a
very promising alternative because the surface wave and current motion
could be handled simultaneously without requiring the computation of
radiation stresses through a wave-averaged model.\\
 It was shown in
\cite{Chen2} that the Green-Naghdi equations (\ref{Green-Naghdi}) are
able to describe partially rotational flows with purely vertical
vorticity. The presence of horizontal vorticity has been considered for
one dimensional surfaces ($d=1$) in \cite{VeeramonySvendsen} for the weakly nonlinear
case and \cite{MusumeciSvendsen} for the fully nonlinear case. These authors
made explicit the contribution of the ``Reynolds'' tensor to the momentum
equation in this case, and showed
 that the momentum equation in (\ref{Green-Naghdi}) must be
modified by the addition of new terms coming from the ``Reynolds''
tensor. The computation of these new terms requires the resolution of
a $1+1$ dimensional transport equation for the vorticity (in
horizontal dimension $d=1$, the stretching term disappears from
(\ref{intreqvort})). The authors derive an approximate explicit solution
to this transport equation under a small amplitude assumption roughly
equivalent to assuming that the vorticity dynamics is weakly
nonlinear.\\
Another recent approach to handle vorticity in shallow water flows has
been proposed in \cite{ZKNPDW}. Their strategy is reminiscent of the
finite element approach: generalizing the approach of
\cite{ShieldsWebster,KBEW} the vertical dependence of the velocity field is projected onto a finite
dimensional basis of functions of $z$ and the ``Reynolds'' tensor ${\bf R}$ is
approximated by its projection onto this basis; computations are
further
simplified by dropping the smallest terms with respect to $\mu$. They
end up with a set of equations on the coordinates of the velocity field in
this basis (these coordinates are
functions of time and of the horizontal variable). Let us also mention
\cite{Constantin} and references therein for the modeling of one
dimensional, periodic or standing rotational water waves.

\medbreak

In this paper, we derive a set of equations, all of them $d$
dimensional, that generalize the Green-Naghdi equations
(\ref{Green-Naghdi}) in presence of vorticity. No assumption is made other
than dropping order $O(\mu^2)$ terms as in the irrotational
theory. Our strategy to get rid of the $d+1$ dimensional dynamics of
the vorticity is inspired by standard turbulence theory. Solving the
vorticity equation is indeed {\it sufficient} to compute the ``Reynolds''
tensor ${\bf R}$, but it is not necessary to do so. One can rather
look for an equation solved by ${\bf R}$. This approach leads to a
cascade of equations (the equation on ${\bf R}$ involves a third order
tensor, which satisfies itself an equation involving a fourth-order
tensor, etc.). The closure of this cascade of equations is one of the
main challenges in turbulence theory; in the present context, we show
that at the order $O(\mu^2)$ of the Green-Naghdi approximation, this
cascade of equations is actually {\it finite}. To be more precise, we
need to introduce first the shear velocity representing the
contribution to the horizontal velocity of the horizontal vorticity,
$$
V_{\rm sh}=\int_z^{\zeta}\bom^\perp_h.
$$
The ``Reynolds'' tensor can then be decomposed into a component (denoted
$E$) due to the
self-interaction of the shear velocity $V_{\rm sh}$, and another one
containing the interaction of the shear velocity with the standard
dispersive vertical dependence of the horizontal velocity due to non-hydrostatic terms,
$$
{\mathbf R}=E+\frac{1}{2}\int_{-1+b}^\zeta\big[ (V^*-\vs)\otimes
(V^*+\vs)+(V^*+\vs)\otimes (V^*-\vs)\big];
$$
with
$$
E=\int_{-1+b}^\zeta \vs\otimes\vs.
$$
These two components are handled separately, and the finite cascade of
equations is derived for $E$; we also show that the other component has a behavior
qualitatively similar to the contribution of the vorticity to non
hydrostatic corrections to the pressure. Three different configurations of
increasing complexity are considered,
\begin{enumerate}
\item {\it One dimensional case with constant vorticity
  $\bom=(0,\omega,0)$}. In this case, the ``Reynolds'' tensor can be
  explicitly computed and the vorticity equation is trivial. Writing
  $\ovV=(\ovv,0)^T$, the
  resulting equations are then\footnote{All the models are derived in
    dimensionless variables in the article, and with an evolution equation for $\ovV$
    rather than $h\ovV$; we chose for the sake of clarity to give the
    dimensional version of these systems in this introduction. }
\begin{equation}\label{GNvortcstintro}
\begin{cases}
\dsp \dt \zeta+\dx (h\ovv)=0,\\
\dsp (1+h {\mathcal T}\frac{1}{h})\big(\dt
(h\ovv)+\dx(h\ovv^2)\big)+gh\dx\zeta +h{\mathcal Q}_1(\ovv)
\\
\hspace{1.5cm}\dsp +\dx\big(\frac{1}{12}h^3\omega^2
\big)+h  {\mathcal C}(\omega h,\ovv) +h {\mathcal C}_b(\omega h,\ovv)=0,
\end{cases}
\end{equation}
where ${\mathcal T}$ and ${\mathcal Q}_1$ are the one-dimensional
versions of the operators defined in (\ref{defTintro}) and
(\ref{defQ1intro}).

    The first term in the third line corresponds to
$\dx E$ and the cascade of equations for $E$ is therefore trivial
(it is equivalent to the mass conservation equation); the last two
terms of the third line gather the non-hydrostatic correction due to
the vorticity and the interaction between the shear velocity and the
dispersive vertical variations (see (\ref{defCCb}) for the definition
of ${\mathcal C}$ and ${\mathcal C}_b$).
\item  {\it One dimensional case with general vorticity}. In this case, the
  ``Reynolds'' tensor cannot be computed explicitly and the cascade of
  equations for $E$ is no longer trivial. For bottom variations of
  medium amplitude (see (\ref{assbottom})), the Green-Naghdi equations
  then become
\begin{equation}\label{GNvort1dfinintro}
\begin{cases}
\dsp \dt \zeta+\dx (h\ovv)=0,\\
\dsp (1\!+\!h{\mathcal T}\frac{1}{h})\big(\dt
(h\ovv)+\dx(h\ovv^2)\big)+gh \dx\zeta \!+\!h {\mathcal Q}_1(\ovv)+\dx E +h{\mathcal C}(v^\sharp,\ovv)=0,
\\
\dsp \dt v^\sharp+\eps \ovv\dx v^\sharp+\eps
v^\sharp\dx \ovv=0,\\
\dsp  \dt E+\eps \ovv\dx E+3\eps E \dx \ovv+\eps\sqrt{\mu}\dx F=0,\\
\dsp \dt F+\eps \ovv\dx F+4\eps F \dx \ovv=0.
\end{cases}
\end{equation}
In these equations, the quantity $F$ that appears in the cascade of
equations for $E$ is the third order self-interaction tensor, while
$v^\sharp$ is introduced to capture corrections due to the non
hydrostatical effects of the vorticity and to the interaction of the shear
velocity with the dispersive vertical variations of the horizontal velocity,
$$
F=\int_{-1+b}^{\zeta}(v^*_{\rm
    sh})^3
\quad
\mbox{ and }\quad
v^\sharp=-\frac{24}{h^3}\int_{-1+
  b}^{\zeta}\int_z^{\zeta}\int_{-1+ b}^z
v^*_{\rm sh}.
$$
\item  {\it Two dimensional case with general vorticity}. This case is
  technically more involved because $E$ is now a $2\times 2$ tensor
  and $F=\int_{-1+b}^\zeta \vs\otimes \vs\otimes \vs$ a $2\times
  2\times 2$ tensor (with coordinates $F_{ijk}$). The main qualitative
  difference is the presence of a source term in the equation for $E$
  that takes into account the interaction between the horizontal and
  vertical components of the vorticity. The equations are
\begin{equation}\label{GNvort1dfinintro2}
\begin{cases}
\dsp \dt \zeta+\nabla\cdot (h\ovV)=0,\\
\dsp (1+h\mathcal{T}\frac{1}{h})\big(\pa_t(h\ovV)
+\nabla\cdot (h\ovV\otimes\ovV)\big)+gh\nabla\zeta+h  {\mathcal
  Q}_1(\ovV)+\nabla\cdot E+h{\mathcal C}(V^\sharp,\ovV)=0,\\
\dsp \pa_t V^\sharp +\ep
(\ovV\cdot\nabla)V^\sharp+\ep(V^\sharp\cdot\nabla)\ovV=0,\\
\dsp \pa_t E + \ep \ovV\cdot\!\nabla \! E +\ep \nabla\!\cdot\!\ovV E +\ep \nabla
\ovV ^T E+\ep E\nabla\ovV+\eps\sqrt{\mu}\nabla\!\cdot\! F
=\eps\sqrt{\mu}
{\mathcal D}(V^\sharp,\ovV) \\
\dsp \pa_t F_{ijk}+\ep \ovV\cdot\nabla F_{ijk}+\ep F_{lkj}\pa_l
\ovV_i+\ep F_{ilk}\pa_l\ovV_j+\ep F_{ijl}\pa_l\ovV_k+\ep \nabla\cdot
\ovV F_{ijk}=0,
\end{cases}
\end{equation}
where the interaction between the horizontal and vertical components
of the vorticity operate through the operator ${\mathcal D}$ (see
(\ref{defD}) for its definition).
\end{enumerate}
A local conservation of energy also holds for these new systems of
equations; we show that
  (\ref{NRJintro}) can be generalized into an equation of the form
\begin{equation}\label{NRJintrorot}
\dt \big({\mathfrak e}+{\mathfrak e}_{\rm rot}\big)+\nabla\cdot \big({\mathfrak
  F}+{\mathfrak F}_{rot}\big)=0,
\end{equation}
(see Remarks
\ref{remNRJcst}, \ref{remNRJgen} and \ref{remNRJ2d} for more
details).

\medbreak

The above Green-Naghdi equations with vorticity allow one to determine
the surface elevation $\zeta$ and the averaged velocity $\ovV$ from
the knowledge of their initial value and of the initial value of the
horizontal vorticity (more precisely, of $V^\sharp$, $E$ and $F$). To
investigate sediment transport for instance, one must be able to
reconstruct the velocity field in the fluid domain. In the
irrotational framework, this structure can be explicitly recovered
from $\zeta$ and $\ovV$. At the level of precision of the Saint-Venant
equations (\ref{StVenant}) (i.e. up to $O(\mu)$ terms in the
approximation for $V$), one gets
\begin{equation}
\label{structvitesseSVirrot}
{\bf V}(t,X,z)=\ovV(t,X),\qquad
w(t,X,z)=0.
\end{equation}
We show that this reconstruction is no longer true in presence of
vorticity and that at the precision of the model, one has for each
level line $\theta$ ($\theta\in [0,1]$),
$$
{\bf V}(t,X,-H_0+b(X)+\theta h(t,X))=\ovV+V_\theta^*,
$$
where $V_\theta^*=V_{\rm sh}^*(t,X,-H_0+b(X)+\theta h(t,X))$ can be
simply determined from its initial value by solving
$$
\dt V_\theta^*+\ovV\cdot \nabla V_\theta^*+V_\theta^*\cdot \nabla\ovV=0.
$$
At the level of precision of the Green-Naghdi equations (i.e. keeping
the $O(\mu)$ terms neglected in the Saint-Venant equations), the
formula for the horizontal velocity becomes
$$
{\bf V}(t,X,-H_0+b(X)+\theta h(t,X))=\ovV+V_\theta^*
+T^*_\theta
\ovV
$$
where the new term $ T^*_\theta$ accounts for dispersive
corrections. These corrections are the same as in the irrotational
theory; the main difference with the first order approximation is that
the equation for $V_\theta^*$ now contains quadratic nonlinear terms
and a source term ${\mathcal S}$, namely,
$$
\dt V^*_\theta+\ovV\cdot \nabla V_\theta^*+V_\theta^*\cdot \nabla
\ovV+V_\theta^*\cdot \nabla
V_\theta^*={\mathcal S}.
$$
The presence of this source term induces an important phenomenon: {\it
  the
creation of a horizontal shear from vertical vorticity}, even if the
initial horizontal vorticity is equal to zero.\\
We finally use these reconstruction formulas for the velocity to
determine the dynamics of the evolution of the vorticity. We show in
particular that the averaged vertical vorticity
$\overline{\bom}_{\mu,v}=\frac{1}{h}\int_{-H_0+b}^\zeta
{\bom}_{\mu,v}$ can be created during the evolution of the
flow throw a {\it mechanism of transfer from horizontal to vertical
  vorticity} which is likely to play an important role for the study
of rip-currents for instance.
\bigbreak

The paper is organized as follows. Section \ref{sect2} is devoted to
an asymptotic analysis of the averaged Euler equations
(\ref{Eulerav}). The first step is to introduce, in \S \ref{sect2.1}, a
dimensionless version of the equations; an asymptotic expansion is
then derived in \S \ref{sect2.2} for the velocity field; this
expansion involves a ``shear velocity'' for which an equation is
derived in \S \ref{sect2.3}, while an expression for the pressure
contribution to (\ref{Eulerav}) is derived in \S \ref{sect2.4}.\\
We then turn to do the derivation of Green-Naghdi type equations in
presence of vorticity; the simplest case of constant vorticity in one
dimension ($d=1$) is first addressed in \S \ref{sectconstant}. The
component $E$ of the ``Reynolds'' tensor ${\bf R}$ and the pressure
contribution can then  be explicitly
computed (see \S \ref{sect3.1} and \S \ref{sect3.2} respectively). The
Green-Naghdi equations with constant vorticity are then derived in \S
\ref{sect3.3}. Section \ref{sect4} then deals with the case of a general
vorticity in dimension $d=1$; there is now a coupling of the momentum
equation with other equations describing the effects of the vorticity
(see \S \ref{sect4.3} to \ref{sect4.5}). The corresponding $1d$ Green-Naghdi
equations with vorticity are derived in \S \ref{sect4.6}. The
two-dimensional case is then handled in Section \ref{sect5}.
We then
explain in Section \ref{sect6} how to reconstruct the velocity field
in the fluid domain and comment on the dynamics of the vorticity. A
first order reconstruction (St-Venant)  is done in \S \ref{sect6.1}
and a second order (Green-Naghdi) one in \S
\ref{sect6.2}. These reconstructions allow us to describe the dynamics
of the vertical vorticity in \S \ref{sect6.3}.
 Finally, a conclusion and perspectives are given
in Section \ref{sect7}.
\section{Asymptotic analysis of the averaged Euler equations (\ref{Eulerav})}\label{sect2}
\subsection{The dimensionless free surface Euler equations}\label{sect2.1}

We non-dimensionize the equations by using several lengths:
the typical amplitude $a_{\rm surf}$ of the waves, the typical
amplitude $a_{\rm bott}$ of the bottom variations, the typical depth  $H_0$, and
the typical horizontal scale $L$. Using these quantities, it is
possible to form three dimensionless parameters,
$$
\eps=\frac{a_{\rm surf}}{H_0},\qquad \beta=\frac{a_{\rm bott}}{H_0},\qquad  \mu=\frac{H_0^2}{L^2};
$$
the parameters $\eps$ and $\beta$ are often called nonlinearity (or
amplitude) and topography
parameters, and the parameter $\mu$ is the {\it shallowness}
parameter.

\begin{remark}\label{remeps}
We are interested here with shallow water flows and
therefore assume that $\mu\ll 1$; on the contrary, we allow for large
amplitude waves and {\it we do not make any
smallness assumption on $\eps$}; it is therefore possible to set $\eps=1$
throughout this article. We however chose to keep track of this
parameter because the simplifications obtained for small amplitude
(weakly nonlinear) or medium amplitude waves can then be performed straightforwardly
(see Remark \ref{remGNCH} for instance). 
\end{remark}

We also use $a_{\rm surf}$, $a_{\rm bott}$, $H_0$ and $L$ to define dimensionless variables and
unknowns (written with a tilde),
$$
\tilde z=\frac{z}{H_0},\qquad \tilde X=\frac{X}{L},\qquad
\tilde\zeta=\frac{\zeta}{a_{\rm surf}}, \qquad \tilde b=\frac{b}{a_{\rm bott}};
$$
the non-dimensionalization of the time variable and of the velocity and pressure
fields is based on the linear analysis of the
equations (see for instance \cite{Lannes_book}, Chapter 1),
$$
\tilde {\bf V}=\frac{{\bf V}}{V_0}, \quad \tilde {\bf w}=\frac{{\bf
    w}}{w_0},\quad
\tilde
t=\frac{t}{t_0},\quad \tilde P=\frac{P}{P_0}
$$
with
$$
V_0=a\sqrt{\frac{g}{H_0}},\quad
w_0=\frac{aL}{H_0}\sqrt{\frac{g}{H_0}},\quad
t_0=\frac{L}{\sqrt{gH_0}},\quad P_0=\rho g H_0.
$$
With these variables and unknowns, and with the notations
$$
\bU^\mu=\left(\begin{array}{c}\sqrt{\mu}{\bf V}\\ {\bf w}\end{array}\right),\qquad
\nabla^\mu=\left(\begin{array}{c}\sqrt{\mu}\nabla\\
    \dz\end{array}\right),
\qquad
N^\mu=\left(\begin{array}{c}-\eps\sqrt{\mu}\nabla\zeta\\ 1\end{array}\right),
$$
and
$$
\curlm=\nabla^\mu\times,\qquad \divem=(\nabla^\mu)^T,\qquad \uU^\mu=(\sqrt{\mu}\uV^T,\uw)^T:=\bU^\mu_{\vert_{z=\eps\zeta}},
$$
the incompressible Euler equations take the form
(omitting the tildes),
\begin{align}
\label{eqUND}
\dt \bU^\mu+\frac{\eps}{\mu} {\bf U}^\mu
\cdot\nabla^\mu\bU^\mu & =-\frac{1}{\eps}\big(\nabla^\mu P + {\bf
  e}_z\big)&\mbox{in }\Omega_t,\\
\nonumber
\divem {\bf U}^\mu&=0&\mbox{in }\Omega_t,
\end{align}
where $\Omega_t$ now stands for the dimensionless fluid domain,
$$
\Omega_t=\{(X,z)\in \R^{d+1},\quad -1+\beta b(X)<z<\eps\zeta(t,X)\}.
$$
Finally, the boundary conditions on the velocity read in dimensionless form,
\begin{eqnarray*}
\dsp \dt \zeta +\nabla\cdot (h\ovV)&=0\quad \mbox{ at the surface,}\\
P&=0\quad \mbox{ at the surface}\\
{\bU}^\mu_{\vert_{z=-1+\beta b}}\cdot N^\mu_b&=0\quad\mbox{ at the bottom},
\end{eqnarray*}
where $N_b^\mu=(-\eps\sqrt{\mu}\nabla b^T,1)^T$ and the dimensionless
versions of $h$ and $\ovV$ are given by
$$
h(t,X):=1+\eps\zeta(t,X)-\beta b(X),\qquad
\ovV(t,X)=\frac{1}{h}\int_{-1+\beta b(X)}^{\eps \zeta(t,X)} {\bf V}(t,X,z)dz.
$$
Decomposing the horizontal velocity into
\begin{equation}
\label{decompVND}
{\bf V}(t,X,z)=\ovV(t,X)+\sqrt{\mu}V^*(t,X,z),
\end{equation}
the dimensionless version of (\ref{Eulerav}) is then given by
\begin{equation}\label{EuleravND}
\begin{cases}
\dsp \dt \zeta+\nabla\cdot (h\ovV)=0,\\
\dsp \dt (h\ovV)+\eps \nabla\cdot (h\ovV\otimes \ovV)+\eps\mu \nabla\cdot (
\int_{-1+\beta b}^{\eps\zeta} V^*\otimes V^*)+\frac{1}{\eps}\int_{-1+\beta b}^{\eps\zeta} \nabla P=0.
\end{cases}
\end{equation}
We therefore need to express the ``rotational Reynolds tensor''
$\int_{-1+\beta b}^{\eps\zeta} V^*\otimes V^*$ and the pressure
contribution $\frac{1}{\eps}\int_{-1+\beta b}^{\eps\zeta} \nabla P$
in terms of $\zeta$ and $\ovV$ in order to obtain a closed set of
equations. This requires a good understanding of the behavior of the
vorticity $\bom_\mu$ defined as
$$
\bom_\mu=\left(
\begin{array}{c}
\frac{1}{\sqrt{\mu}}(\dz {\bf V}^\perp-\nabla^\perp {\bf w})\\
-\nabla\cdot {\bf V}^\perp
\end{array}\right)=\frac{1}{\mu} \curlm {\bf U}^\mu.
$$
We shall assume here that $\bom_\mu$ is of order $O(1)$ with respect
to $\mu$ (in the terminology of \cite{Teshukov,RG1}, we consider therefore
weakly sheared flows). It is rigorously shown in \cite{CastroLannes}
that in this regime, $\zeta$, $\ovV$ and $\bom_\mu$ remain $O(1)$
quantities during the time evolution of the flow. All the formal
asymptotic descriptions made throughout this article have therefore a
firm basis, and the models derived here could be rigorously justified
using the procedure hinted in \cite{CastroLannes} (and carried through
for the Saint-Venant equations with vorticity).

\subsection{Asymptotic expansion of the velocity field}\label{sect2.2}

As explained above, the quantities $\zeta$, $\ovV$ and $\bom_\mu$ are
all of order $O(1)$ with respect to $\mu$. An asymptotic description
of $\bU^\mu$ can then be found by considering the boundary value
problem
\begin{equation}\label{divrotND}
\left\lbrace
\begin{array}{lll}
\curlm \bU^\mu &= \mu\bom_\mu&\mbox{ in }\quad \Omega\\
\divem \bU^\mu &=0&\mbox{ in }\quad \Omega\\
%\Vp^\mu&= \nabla\psi+\nabla^\perp \Delta^{-1}(\uom_\mu\cdot N^\mu)&\mbox{ at the surface}\\
U^\mu_b\cdot N^\mu_b&=0 &\mbox{ at the bottom}
\end{array}\right.
\end{equation}
(the subscript $b$ is used to denote quantities evaluated at the
bottom). 
\begin{remark}Only the ``rotational'' part of $\bU^\mu$ is fully determined
  from this boundary value problem. The ``irrotational part'' of
  $\bU^\mu$ is determined from the tangential component of the
  velocity at the interface; more precisely, one can show that
$\bU^\mu_\parallel:={\bf V}_{\vert_{z=\eps\zeta}}+\eps {\bf
  w}_{\vert_{z=\eps\zeta}}\nabla\zeta$ can be written under the form
$$
\bU^\mu_\parallel=\nabla\psi
+\nabla^\perp\Delta^{-1}({\bom_{\mu}}_{\vert_{z=\eps\zeta}}\cdot N^\mu),
$$
for some scalar function $\psi$ defined over $\R^2$. This function
fully determines the ``irrotational'' part of $\bU^\mu$ which is given
by $\nabla^\mu \Phi$, with
$$
\begin{cases}
(\dz^2+\mu \Delta)\Phi=0\qquad \mbox{ \rm in }\quad \Omega\\
\Phi_\surf=\psi,\qquad \dz\Phi_\bott=0
\end{cases}
$$
(and which obviously leaves (\ref{divrotND}) unchanged).
We refer to \cite{CastroLannes} for more details; it is in particular
shown in this reference that $\zeta$, $\psi$ and $\bom_\mu$ remain
uniformly bounded with respect to $\mu$ during the time evolution of
the flow. The analysis of the full boundary value problem then shows that ${\bf V}$, ${\bf
  w}$, $\ovV$ etc. also remain bounded, which justifies the formal
asymptotics made here.
\end{remark}

We show here how to construct an approximate solution to this system
of equations. Replacing
\begin{equation}\label{ansU}
\bU^\mu=\left(\begin{array}{c} \sqrt{\mu}{\bf V} \\ {\bf
      w}\end{array}\right)=\left(\begin{array}{c} \sqrt{\mu}\ovV+\mu V^* \\ \mu\tilde{\bf w}\end{array}\right)
\end{equation}
in the equations (\ref{divrotND}) and using the fact that $\dz
\ovV=0$, we obtain
\begin{equation}\label{divrotND1}
\left\lbrace
\begin{array}{lll}
\dz V^*-\sqrt{\mu}\nabla \tilde {\bf w} &= -\bom_{\mu,h}^\perp&
\mbox{ in }\quad \Omega\\
\nabla^\perp\cdot \ovV+\sqrt{\mu}\nabla^\perp\cdot V^*&=\bom_{\mu,v} &
\mbox{ in }\quad \Omega\\
\nabla\cdot \ovV+\sqrt{\mu}\nabla\cdot V^*+\dz \tilde {\bf w} &=0&\mbox{ in }\quad \Omega\\
%\Vp^\mu&= \nabla\psi+\nabla^\perp \Delta^{-1}(\uom_\mu\cdot N^\mu)&\mbox{ at the surface}\\
\tilde{\bf w}_b-\beta \nabla b\cdot (\ovV+\sqrt{\mu}V^*_b)&=0 &\mbox{ at the bottom},
\end{array}\right.
\end{equation}
where $\bom_{\mu,h}$ and $\bom_{\mu,v}$ denote respectively the
horizontal and vertical components of the vorticity $\bom_\mu$. From
the third and last equations of (\ref{divrotND1}), we first get an
expression for $\tilde{\bf w}$ in terms of $\ovV$ and $V^*$,
\begin{equation}\label{exprw}
\tilde{\bf w}=-\nabla\cdot \big[(1+z-\beta b)\ovV
\big]-\sqrt{\mu}\nabla\cdot \int_{-1+\beta b}^{z}V^*.
\end{equation}
Replacing this expression in the first equation of (\ref{divrotND1})
then gives
$$
\dz V^*=\sqrt{\mu}\nabla \tilde {\bf w}-\bom_{\mu,h}^\perp
$$
and therefore
$$
V^*=\sqrt{\mu}\Big(\int_{z}^{\eps\zeta}\nabla\nabla\cdot \big[(1+z'-\beta b)\ovV
\big]dz'\Big)^*+\mu \Big(\int_{z}^{\eps\zeta}\nabla\nabla\cdot
\int_{-1+\beta b}^zV^*\Big)^*+\Big(\int_z^{\eps\zeta}\bom_{\mu,h}^\perp\Big)^*.
$$
Defining the operators $T=T[\beta
b,\eps\zeta]$ and $T^*=T^*[\beta
b,\eps\zeta]$ by
\begin{eqnarray*}
T[\beta
b,\eps\zeta] W=\int_{z}^{\eps\zeta}\nabla\nabla\cdot
\int_{-1+\beta b}^{z'}W,\qquad \mbox{ and }\qquad
T^*[\beta
b,\eps\zeta] W=(T[\beta
b,\eps\zeta] W)^*,
\end{eqnarray*}
and the ``shear'' velocity
$V_{\rm sh}$ by
\begin{equation}\label{defshear}
V_{\rm sh}=\int_z^{\eps\zeta}\bom_{\mu,h}^\perp,
\end{equation}
we can rewrite the above identity under the form
$$
(1-\mu T^*)V^*=\sqrt{\mu}T^*\ovV+V^*_{\rm sh}.
$$
Note that this is an {\it exact} identity. Applying $(1+\mu T^*)$ on
both sides, we obtain the following approximation of order $O(\mu^{3/2})$,
\begin{equation}\label{approxV*}
V^*=\sqrt{\mu}T^*\ovV +(1+\mu T^*)V^*_{\rm sh}+O(\mu^{3/2}).
\end{equation}
Together with (\ref{decompVND}), this yields the following $O(\mu)$
and $O(\mu^2)$ approximations
for the horizontal velocity field ${\bf V}$,
\begin{align}
\label{Vmu}
{\bf V}&=\ovV+\sqrt{\mu}V^*_{\rm sh}+O(\mu),\\
\label{Vmu2}
{\bf V}&=\ovV+\sqrt{\mu}V^*_{\rm sh}    +\mu T^*\ovV +\mu^{3/2} T^*V^*_{\rm sh}+O(\mu^2).
\end{align}
\begin{remark}
Remark that even at order $O(\mu)$ the flow cannot be assumed to be
columnar (i.e. that its horizontal velocity is independent of $z$),
which is in sharp contrast with the rotational case.
\end{remark}
\begin{remark}
Since $\ovV$ does not depend on $z$, one can compute explicitly
$$
T^*\ovV=-\frac{1}{2}\big((z+1-\beta
b)^2-\frac{h^2}{3}\big)\nabla\nabla\cdot \ovV+\beta
\big(z-\eps\zeta+\frac{1}{2}h\big)\big[\nabla b\cdot
\nabla\ovV+\nabla(\nabla b\cdot \ovV)\big].
$$
\end{remark}
Finally, for the vertical component ${\bf w}$ of the velocity, we use (\ref{approxV*}) together with (\ref{ansU}) and (\ref{exprw}) to get
\begin{align}
\label{wmu}
{\bf w}&=O(\mu),\\
\label{wmu2}
{\bf w}&=-\mu\nabla\cdot \big[(1+z-\beta b)\ovV
\big]-\mu^{3/2}\nabla\cdot \int_{-1+\beta b}^{z}V^*_{\rm sh}+O(\mu^2).
\end{align}

\subsection{An equation for the shear velocity}\label{sect2.3}

The purpose of this section is to derive an approximate equation solved by
the shear velocity $V^*_{\rm sh}$ defined in (\ref{defshear}). We recall that the dimensionless vorticity equations is given by
$$
\dt \bom_\mu+\frac{\eps}{\mu}\bU^\mu\cdot \nabla^\mu
\bom_\mu=\frac{\eps}{\mu}\bom_\mu\cdot \nabla^\mu \bU^\mu,
$$
so that its horizontal component $\bom_{\mu,h}$ solves
\begin{equation}\label{eqrotNDh}
\dt \bom_{\mu,h}+{\eps}{\bf V}\cdot \nabla
\bom_{\mu,h}+\frac{\eps}{\mu}{\bf w}\dz
{\bom_{\mu,h}}=\eps\bom_{\mu,h}\cdot\nabla {\bf
  V}+\frac{\eps}{\sqrt{\mu}}{\bom_{\mu,v}}\dz{\bf V}.
\end{equation}
Using (\ref{Vmu}) and (\ref{wmu2}), this yields
\begin{align}
\nonumber
\dt \bom_{\mu,h}+{\eps}{\ovV}\cdot \nabla
\bom_{\mu,h}&-\eps \nabla\cdot [(1+z-\beta b)\ovV]\dz
{\bom_{\mu,h}}\\
\label{ord1}
&=\eps\bom_{\mu,h}\cdot\nabla
{\ovV}-\eps({\nabla^\perp\cdot \ovV})\bom_{\mu,h}^\perp+O(\eps\sqrt{\mu})
\end{align}
(as explained in Remark \ref{remeps}, we keep track of the dependence
on $\eps$, but no smallness assumption is made on this parameter, and
one can set $\eps=1$ everywhere). Recalling that $V_{\rm sh}=\int_z^{\eps\zeta}\bom_{\mu,h}^\perp$, we
can integrate this equation with respect to $z$ to get
\begin{align*}
\nonumber
\dt V_{\rm sh}&+\eps\ovV\cdot \nabla V_{\rm sh}+\eps(\nabla\cdot
\ovV)V_{\rm sh}-(\dt \zeta+\eps\ovV\cdot \nabla
\zeta){\bom_{\mu,h}^\perp}_\surf\\
&-\eps \Big[\nabla\cdot [(1+z-\beta b)\ovV]
{\bom_{\mu,h}^\perp}\Big]^{\eps\zeta}_{z} =-\eps V_{\rm sh}^\perp\cdot \nabla\ovV^\perp-\eps
(\nabla^\perp\cdot \ovV)V_{\rm sh}^\perp+O(\eps\sqrt{\mu}).
\end{align*}
Using the fact that $\dt \zeta+\nabla\cdot (h\ovV)=0$ we therefore get
\begin{align*}
\dt V_{\rm sh}+\eps\ovV\cdot \nabla V_{\rm sh}+\eps(\nabla\cdot
\ovV)V_{\rm sh}&+\eps \nabla\cdot [(1+z-\beta b)\ovV]
{\bom_{\mu,h}^\perp} \\
&=-\eps V_{\rm sh}^\perp\cdot \nabla\ovV^\perp-\eps
(\nabla^\perp\cdot \ovV)V_{\rm sh}^\perp+O(\eps\sqrt{\mu});
\end{align*}
recalling the vectorial identity
\begin{equation}\label{vectorialid}
(\nabla\cdot A)B+B^\perp\cdot \nabla A^\perp+(\nabla^\perp\cdot
A)B^\perp=B\cdot \nabla A,
\end{equation}
this equation can be rewritten
\begin{eqnarray}
\label{ord2}
\dt V_{\rm sh}+\eps\ovV\cdot \nabla V_{\rm sh}+\eps V_{\rm sh}\cdot
\nabla \ovV-\eps \nabla\cdot [(1+z-\beta b)\ovV]
\dz V_{\rm sh} =O(\eps\sqrt{\mu}).
\end{eqnarray}
Integrating this equation, one readily obtains the following equation
for the average shear velocity $\ovV_{\rm sh}$,
$$
\dt \ovV_{\rm sh}+\eps \ovV\cdot \nabla \ovV_{\rm sh}+\eps \ovV_{\rm
  sh}\cdot \nabla \ovV=O(\eps\sqrt{\mu}).
$$
Subtracting this equation to (\ref{ord2}), we finally get the
following evolution equation on $V_{\rm sh}^*$,
\begin{equation}\label{eqVsh}
\dt V_{\rm sh}^*+\eps\ovV\cdot \nabla V_{\rm sh}^*+\eps V_{\rm sh}^*\cdot
\nabla \ovV-\eps \nabla\cdot [(1+z-\beta b)\ovV]
\dz V_{\rm sh}^* =O(\eps\sqrt{\mu}).
\end{equation}
Note that we shall derive later a more precise evolution equation for
(\ref{eqVsh}) by making explicit the $O(\eps\sqrt{\mu})$ terms (see \S
\ref{secteqvsharp2} below).

\subsection{Asymptotic expansion of the pressure field}\label{sect2.4}

The vertical component of the Euler equation (\ref{eqUND}) is given
by
$$
\dt {\bf w}+\eps {\bf V}\cdot \nabla {\bf w}+\frac{\eps}{\mu}{\bf
  w}\dz {\bf w}=-\frac{1}{\eps}(\dz P+1);
$$
since moreover $P$ vanishes at the surface, we get that
\begin{align*}
\nonumber
\frac{1}{\eps}\nabla P&=\nabla\int_{z}^{\eps \zeta}
(-\frac{1}{\eps}\dz P)\\
&=\nabla\zeta+\nabla\int_{z}^{\eps\zeta}\Big(\dt {\bf w}+\eps {\bf V}\cdot \nabla {\bf w}+\frac{\eps}{\mu}{\bf
  w}\dz {\bf w}\Big)
\end{align*}
and therefore
\begin{equation}
\label{approxP}
\frac{1}{\eps}\int_{-1+\beta b}^{\eps\zeta}\nabla
P=h\nabla\zeta+\int_{-1+\beta b}^{\eps\zeta}\nabla\int_{z}^{\eps\zeta}\Big(\dt {\bf w}+\eps {\bf V}\cdot \nabla {\bf w}+\frac{\eps}{\mu}{\bf
  w}\dz {\bf w}\Big).
\end{equation}
The first term in the right-hand-side corresponds to the hydrostatic
pressure. We still need an expansion of the non-hydrostatic terms with
respect to $\mu$; for the sake of clarity, this computation is
performed in \S \ref{sect3.2} for the one dimensional case with
constant vorticity, in \S \ref{sectcontribP1d} for the one dimensional case
with general vorticity, and in \S \ref{sect5.2} for the two
dimensional case.

\begin{remark}\label{remSV}
A direct consequence of (\ref{approxP}) and (\ref{wmu}) is that
$$
\frac{1}{\eps}\int_{-1+\beta b}^{\eps\zeta}\nabla
P=h\nabla\zeta+O(\mu).
$$
An order $O(\mu)$ approximation of (\ref{EuleravND}) is therefore
provided by the Saint-Venant (or nonlinear shallow water) equations
\begin{equation}\label{StV}
\begin{cases}
\dsp \dt \zeta+\nabla\cdot (h\ovV)=0,\\
\dsp \dt (h\ovV)+h\nabla\zeta+\eps \nabla\cdot (h\ovV\otimes \ovV)=0;
\end{cases}
\end{equation}
these equations are exactly the same as in the irrotationnal
setting. As we shall see in the next sections, the rotational terms
affect the $O(\mu)$ terms in (\ref{EuleravND}) and, consequently, the
Green-Naghdi equations differ from the standard irrotational version
when vorticity is present.
\end{remark}

\section{The $1d$ Green-Naghdi equations with constant vorticity} \label{sectconstant}

In dimension $d=1$, one can consider flows with constant vorticity,
\begin{equation}\label{not1d}
\bU^\mu=\left(\begin{array}{c}  \sqrt{\mu} v \\ 0 \\ {\bf
      w}  \end{array}\right),
\qquad
\bom_\mu=\left(\begin{array}{c}  0 \\ \omega \\0 \end{array}\right),
\end{equation}
(so
that $-\bom_{\mu,h}^\perp=(\omega,0)^T$), and with $\omega=\dz
v^*-\sqrt{\mu}\dx \tilde{\bf w}\equiv \mbox{constant}$. The fact that
the vorticity remains constant implies that the wave does not affect
the underlying current; however, as we show in this section, the
current modifies the motion of the waves.

Since $\omega$ is constant, one deduces from (\ref{defshear}) that $V^*_{\rm sh}=(v^*_{\rm sh},0)^T$ with
\begin{equation}\label{exprVstar}
v^*_{\rm sh}=-\big((\eps\zeta-z)-\frac{1}{2}h\big)\omega
\end{equation}
and, writing $\ovV=(\ovv,0)^T$, (\ref{wmu2}) takes therefore the form
\begin{equation}\label{wrotcst}
{\bf w}=-\mu\dx\big( (1+z-\beta b)
\ovv\big)+\mu^{3/2}\frac{1}{2}\omega\dx\big( (\eps\zeta-z)(z+1-\beta b)\big).
\end{equation}

\subsection{Computation of the ``rotational Reynolds tensor''
  contribution}\label{sect3.1}

Using (\ref{exprVstar}) and (\ref{wrotcst}),  the contribution of the ``rotational Reynolds-tensor''
to (\ref{EuleravND}) can be written (dropping $O(\mu^2)$ terms)
\begin{align}
\nonumber
\eps\mu\dx
\int_{-1+\beta b}^{\eps\zeta} \abs{v^*}^2=&\eps\mu\dx
\int_{-1+\beta b}^{\eps\zeta} \abs{v^*_{\rm sh}}^2+2\eps\mu^{3/2}\dx
\int_{-1+\beta b}^{\eps\zeta}
v^*_{\rm sh} T\ovv\\
\nonumber
=&\eps\mu\omega^2\dx
\int_{-1+\beta b}^{\eps\zeta}
\big((\eps\zeta-z)-\frac{1}{2}h\big)^2\\
\nonumber
&-2\eps\mu^{3/2}\omega\dx
\int_{-1+\beta b}^{\eps\zeta}
\big((\eps\zeta-z)-\frac{1}{2}h\big)
\int_{z}^{\eps\zeta}\dx^2\big((1+z'-\beta b)\ovV\big)\\
\label{contribRcst}
=&\frac{\eps\mu}{12}\omega^2\dx (h^3)
-\frac{\eps\mu^{3/2}}{12}\omega\dx \big[h^3\big(h\dx^2
\ovv\!-\!2\beta\ovv\dx^2b\!-\!4\beta\dx\ovv\dx b\big)\big].
\end{align}
Note that this expression only depends on $\zeta$ and $\ovv$ (and on the
constant vorticity $\omega$ and the bottom parametrization $b$).
\subsection{Computation of the pressure
  contribution}\label{sect3.2}

Similarly, we can use (\ref{exprVstar}) and (\ref{wrotcst}) to write
the pressure contribution to
(\ref{approxP})  under the form (dropping again the $O(\mu^2)$ terms)
\begin{align}
\nonumber
\frac{1}{\eps}\!\int_{-1+\beta b}^{\eps\zeta}\hspace{-3mm}\dx
P=&h\dx\zeta+\int_{-1+\beta b}^{\eps\zeta}\dx\int_{z}^{\eps\zeta}\Big(\dt {\bf w}+\eps v\dx {\bf w}+\frac{\eps}{\mu}{\bf
  w}\dz {\bf w}\Big)\\
\nonumber
=&h\dx\zeta+\mu h {\mathcal T} \big((\dt \ovv+\eps
\ovv\dx\ovv\big)+\mu\eps h{\mathcal Q}_1(\ovv)\\
\nonumber
&-\eps\mu^{3/2}\frac{1}{2}\omega h {\mathcal T}\big(\beta\ovv\dx b-h\dx v\big)-\eps\mu^{3/2}\frac{1}{12}\omega\dx
(h^4\dx^2\ovv+6h^3\dx\zeta\dx\ovv)\\
\label{contribPcst}
&+\!\beta\mu^{3/2}\frac{1}{12}\big[\omega h^2\dx \ovv\big(h\dx^2b\!-\!6\dx\zeta\dx
b\big)+\!3\beta h \ovv \dx b(h\dx^2 b+2\dx b \dx \zeta)\big],
\end{align}
where ${\mathcal T}$ and ${\mathcal Q}_1$ are defined in (\ref{defTintro}) and (\ref{defQ1intro}), and where we
used the fact that $\dt \zeta+\ovv\dx \zeta= \beta \ovv\dx b-h\dx v$.
As for the Reynolds-tensor contribution (\ref{contribRcst}), this
expression only depends on the variables $\zeta$ and $\ovV$.
\begin{remark}
 We actually work here with the dimensionless form of the operators ${\mathcal T}$ and ${\mathcal Q}_1$, formally
  obtained by replacing $\zeta$ by $\eps\zeta$, $b$ by $\beta b$ and
  with $h=1+\eps\zeta-\beta b$. In dimension $d=1$, this gives
\begin{align*}
{\mathcal T}v&=-\frac{1}{3h}\dx(h^3\dx
v)+\beta \frac{1}{2h}\big[\dx (h^2\dx b v)-h^2\dx b\dx
v\big]+\beta^2 (\dx b)^2 v\\
{\mathcal Q}_1(v)&=\frac{2}{3h}\dx \big(h^3 (\dx v)^2\big)+\beta h (\dx
v)^2\dx b+\frac{\beta}{2h}\dx \big(h^2 v^2 \dx^2 b\big)+\beta ^2v^2\dx b\dx^2b.
\end{align*}
\end{remark}

\subsection{The Green-Naghdi model}\label{sect3.3}

Gathering (\ref{contribRcst}) and (\ref{contribPcst}), we obtain
\begin{align*}
\mu\eps\dx
\int_{-1}^{\eps\zeta} \abs{v^*}^2&+\frac{1}{\eps}\int_{-1}^{\eps\zeta}\dx
P=h\dx\zeta\\
&+\mu h{\mathcal T} (\dt
\ovv+\eps\ovv\dx\ovv)+\eps\mu h{\mathcal Q}_1(\ovv)
+\eps\mu\frac{1}{12}\omega^2\dx (h^3)\\
&-\eps\mu^{3/2}\frac{1}{6}\dx\big(2\omega h^4\dx^2\ovv+\dx(\omega
h^4)\dx\ovv\big)
+\eps\beta\mu^{3/2}\frac{\omega}{3}\big[\dx\big(h^3\dx^2 b\ovv\big)\!+\!h^3\dx^2b\dx \ovv\big]
\end{align*}
(up to $O(\mu^2)$ terms); plugging this approximation into
(\ref{EuleravND}), we obtain the Green-Naghdi equations in dimension
$d=1$ and with constant vorticity,
\begin{equation}\label{GNvortcst}
\begin{cases}
\dsp \dt \zeta+\dx (h\ovv)=0,\\
\dsp (1+\mu{\mathcal T})\big(\dt
\ovv+\eps\ovv\dx\ovv\big)+\dx\zeta +\eps\mu{\mathcal Q}_1(\ovv)
\\
\hspace{0.4cm}\dsp +\eps\mu\frac{1}{h}\dx\big(\frac{1}{12}h^3\omega^2
\big)+\eps\mu^{3/2}  {\mathcal C}(\omega h,\ovv) +\eps\beta \mu^{3/2}  {\mathcal C}_b(\omega h,\ovv)=0,
\end{cases}
\end{equation}
where ${\mathcal T}$ and ${\mathcal Q}_1$ are as in (\ref{defTintro}) and (\ref{defQ1intro}), and with
\begin{equation}\label{defCCb}
\forall v^\sharp, \forall \ovv,
\qquad
\begin{array}{lcl}
\dsp {\mathcal C}(v^\sharp,\ovv)
&=&\dsp -\frac{1}{6h}\dx\big(2 h^3 v^\sharp\dx^2\ovv+\dx( h^3
v^\sharp)\dx\ovv\big)\vspace{.5mm}\\
\dsp  {\mathcal C}_b(v^\sharp,\ovv)
&=&\dsp \frac{1}{3h}\big[\dx\big(h^2 v^\sharp\dx^2
b\ovv\big)+h^2v^\sharp\dx^2b\dx \ovv\big].
\end{array}
\end{equation}
\begin{remark}\label{remfluxC}
One can notice that
$$
h{\mathcal C}(v^\sharp,\ovv)\cdot \ovv=\dx
{\mathfrak F}_{{\mathcal C}}
\quad \mbox{ and }\quad
h{\mathcal C}_b(v^\sharp,\ovv)\cdot\ovv=\dx
{\mathfrak F}_{{\mathcal C}_b},
$$
where the fluxes ${\mathfrak F}_{{\mathcal C}}$ and ${\mathfrak
  F}_{{\mathcal C}_b}$ are given by
$$
{\mathfrak F}_{{\mathcal C}}=\frac{h^3}{6}\big(
(\dx\ovv)^2-\ovv\dx^2 v\big)v^\sharp-\frac{1}{6}\dx (h^3 v^\sharp \dx
v)\ovv
\quad\mbox{ and }\quad
{\mathfrak F}_{{\mathcal C}_b}=h^2 v^\sharp \dx^2 b \ovv^2.
$$
\end{remark}
\begin{remark}\label{remNRJcst}
Remark \ref{remfluxC} can be used to derive a local equation for the
conservation of energy; one finds indeed that (\ref{NRJintro}) can be
generalized in presence of a constant vorticity into
\begin{equation}\label{NRJconstvort}
\dt \big({\mathfrak e}+{\mathfrak e}_{\rm rot}\big)+\dx
\big({\mathfrak F}+{\mathfrak F}_{\rm rot}\big)=0,
\end{equation}
with (in dimensional form)
$$
{\mathfrak e}_{\rm rot}=  \frac{1}{24}\omega^2 h^3
\quad\mbox{ and }\quad
{\mathfrak F}_{\rm rot}=\frac{1}{8}\omega h^3+{\mathfrak F}_{{\mathcal C}}+{\mathfrak F}_{{\mathcal C}_b}.
$$
\end{remark}
\section{The $1d$ Green-Naghdi equations with general vorticity}\label{sect4}

In dimension $d=1$ with a non constant vorticity, we still use the
notations (\ref{not1d}), but $\omega=\omega(t,x,z)$ now depends on
space and time. Contrary to what happens in the constant vorticity
case studied in Section \ref{sectconstant}, we show here that there
is a nontrivial wave-current interaction in the sense that the
underlying current is now affected by the motion of the waves.

Still denoting $V^*_{\rm sh}=(v^*_{\rm sh},0)^T$ , one deduces from (\ref{defshear}) that 
\begin{equation}\label{exprVstar1d}
v^*_{\rm sh}=-\Big(
\int_z^{\eps\zeta}\omega\big)^*.
\end{equation}
As in Section \ref{sectconstant} we compute the contribution to
(\ref{EuleravND}) of the ``rotational Reynolds-tensor'' and of the
pressure; the difference with the case of constant vorticity addressed
in Section \ref{sectconstant} is that the component $E$ of the ``Reynolds'' tensor
cannot be computed explicitly; we show that it can be determined
through the resolution of a finite cascade of equations.

In order to simplify the computations, we shall make from now on the
following assumption on the amplitude of the bottom variations
\begin{equation}\label{assbottom}
\mbox{Medium amplitude bottom variations:} \qquad \beta=O(\sqrt{\mu}).
\end{equation}
This assumption leads to simpler models without neglecting any essential
mechanisms of wave-current interaction; for instance, in the particular case of
constant vorticity, this assumption allows one to neglect the term
$\eps\beta\mu^{3/2}{\mathcal C}_b(\omega h,\ovv)$ in (\ref{GNvortcst}).

\subsection{Computation of the ``rotational Reynolds tensor''
  contribution}\label{sect4.1}

Proceeding as for (\ref{contribRcst}), we write (dropping $O(\mu^2)$ terms)
\begin{align*}
\eps\mu\dx
\int_{-1+\beta b}^{\eps\zeta} \abs{v^*}^2&=\eps\mu\dx
\int_{-1+\beta b}^{\eps\zeta} \abs{v^*_{\rm sh}}^2+2\eps\mu^{3/2}\dx
\int_{-1+\beta b}^{\eps\zeta}
v^*_{\rm sh} T\ovv\\
&=\eps\mu\dx
\int_{-1+\beta b}^{\eps\zeta} \abs{v^*_{\rm sh}}^2+2\eps\mu^{3/2}\dx
\int_{-1+\beta b}^{\eps\zeta}
v^*_{\rm sh} \int_{z}^{\eps\zeta}\dx^2\big((1+z'-\beta b)\ovv\big)\\
&=\eps\mu\dx
\int_{-1+\beta b}^{\eps\zeta} \abs{v^*_{\rm sh}}^2-\eps\mu^{3/2}\dx
\int_{-1+\beta b}^{\eps\zeta}
\int_{z}^{\eps\zeta}v^*_{\rm sh} \dx^2\big((1+z'-\beta b)\ovv\big),
\end{align*}
the last line stemming from the identity
$$
2\int_{-1+\beta b}^\zeta W^*\int_z^{\eps\zeta}\dx^2\big((1+z'-\beta
b)\overline{W}\big)=-\int_{-1+\beta
  b}^{\eps\zeta}\int_{z}^{\eps\zeta}W^*\dx^2\big((1+z'-\beta b)\overline{W}\big).
$$
Introducing
$$
v^\sharp:=-\frac{24}{h^3}\int_{-1+\beta
  b}^{\eps\zeta}\int_z^{\eps\zeta}\int_{-1+\beta b}^z
v^*_{\rm sh},
$$
one computes that
\begin{equation}\label{idvs1}
\int_{-1+\beta b}^{\eps\zeta}
\int_{z}^{\eps\zeta}v^*_{\rm sh} \dx^2\big((1+z'-\beta b)\ovv\big)
=\frac{h^3}{12}v^\sharp\dx^2 \ovv+O(\beta)
\end{equation}
so that under the assumption
(\ref{assbottom}), one readily deduces that
\begin{eqnarray}
\eps\mu\dx
\int_{-1+\beta b}^{\eps\zeta} \abs{v^*}^2&
\label{contribR1d}
=\eps\mu\dx
\int_{-1+\beta b}^{\eps\zeta} \abs{v^*_{\rm sh}}^2-\eps\mu^{3/2} \frac{1}{12}\dx(h^3v^\sharp\dx^2\ovv)+O(\mu^2).
\end{eqnarray}
\subsection{Computation of the pressure
  contribution}\label{sectcontribP1d}

As for (\ref{contribPcst}), we write
\begin{align}
\nonumber
\frac{1}{\eps}\int_{-1+\beta b}^{\eps\zeta}\dx
P\dsp &=h\dx\zeta+\int_{-1+\beta b}^{\eps\zeta}\dx\int_{z}^{\eps\zeta}\Big(\dt {\bf w}+\eps v\dx {\bf w}+\frac{\eps}{\mu}{\bf
  w}\dz {\bf w}\Big)\\
\label{presscontr}
&=h\dx\zeta+ \mu A_{(1)}+\mu^{3/2}A_{(3/2)}+O(\mu^{2}),
\end{align}
with
\begin{equation}\label{presscontr1}
A_{(1)}= h {\mathcal T} \big((\dt \ovv+\eps
\ovv\dx\ovv\big)+\eps h{\mathcal Q}_1(\ovv).
\end{equation}
For $A_{(3/2)}$, it is convenient to introduce the operator $\widetilde{\mathcal T}=\widetilde{\mathcal T}[\eps\zeta,\beta b]$
defined\footnote{We use here the one dimensional version $
\widetilde{\mathcal T}W=-\frac{1}{h}\int_{-1+\beta
    b}^{\eps\zeta}\dx\int_{z}^{\eps\zeta}\dx
  \int_{-1+\beta b}^z W.
$} as
$$
\widetilde{\mathcal T}[\eps\zeta,\beta b]W=-\frac{1}{h}\int_{-1+\beta
    b}^{\eps\zeta}\nabla\int_{z}^{\eps\zeta}\nabla\cdot
  \int_{-1+\beta b}^z W
$$
(in particular, if $W=\overline{W}$ does not depend on $z$, then
$\widetilde{T}\overline{W}={\mathcal T}\overline{W}$). 

We then have
\begin{align*}
A_{(3/2)}=&h\widetilde{\mathcal T}\dt v^*_{\rm
  sh}\\
&+\eps\int_{-1+\beta
  b}^{\eps\zeta}\dx\int_{z}^{\eps\zeta}\big(\ovv\dx(-\dx\int_{-1+\beta
b}^zv^*_{\rm
sh})+v^*_{\rm sh}\dx (-\dx((1+z'-\beta b)\ovv))\big)\\
&+\eps\int_{-1+\beta
  b}^{\eps\zeta}\dx\int_{z}^{\eps\zeta}\big(\!\!-\dx((1+z'-\beta b)\ovv)(-\dx
v^*_{\rm sh})-\int_{-1+\beta b}^z\dx v^*_{\rm sh}(-\dx \ovv)\big),
\end{align*}
up to $O(\sqrt{\mu})$ terms.
We can notice that the contribution of the bottom of  order
$O(\beta)$ in $A_{(3/2)}$ and therefore of order $O(\mu^{3/2}\beta)$
in (\ref{presscontr}). For bottoms of medium amplitude as
in (\ref{assbottom}), this contribution is of order $O(\mu^2)$ and can
therefore be neglected at the precision of the model. We therefore
take $b=0$ in the following computations, which allows us to write
(dropping $O(\sqrt{\mu},\beta)$ terms)
\begin{align}
\nonumber
A_{(3/2)}=&h\widetilde{\mathcal T}(\dt v^*_{\rm
  sh}+\eps\ovv\dx v^*_{\rm sh})
+2\eps\int_{-1}^{\eps\zeta}\dx\int_{z}^{\eps\zeta}(\dx
\ovv)\dx\int_{-1}^{z'}v^*_{\rm
sh}\\
\nonumber
& -\eps\int_{-1}^{\eps\zeta}\dx \int_z^{\eps\zeta} v^*_{\rm sh}
\dx^2 ((1+z')\ovv)\\
\label{presscontr2}
&+\eps\int_{-1}^{\eps\zeta}\dx \int_z^{\eps\zeta} (\dx v^*_{\rm sh})
\dx ((1+z')\ovv).
\end{align}
In order to get a simpler expression for $B_{(3/2)}$, let us apply
$\widetilde{\mathcal T}$ to (\ref{eqVsh}) to get, up to
$O(\sqrt{\mu})$ terms,
\begin{align*}
h\widetilde{\mathcal T}\big(\dt v_{\rm sh}^*&+\eps\ovv\dx v_{\rm sh}^*+\eps v_{\rm sh}^*\dx \ovv\big)=-\eps\int_{-1}^\zeta\dx \int_{z}^{\eps\zeta}\dx \int_{-1}^{z'} \dx [(1+z)\ovv]
\dz v_{\rm sh}^* \\
&=\eps\int_{-1}^\zeta\dx \int_{z}^{\eps\zeta}\dx \int_{-1}^{z'} \dx \ovv
 v_{\rm sh}^*-\eps\int_{-1}^\zeta\dx \int_{z}^{\eps\zeta}\dx [\dx ((1+z')\ovv)
 v_{\rm sh}^*]
\end{align*}
and therefore
\begin{align*}
h\widetilde{\mathcal T}\big(\dt v_{\rm sh}^*+\eps\ovv\dx v_{\rm sh}^*\big)=&-2h\widetilde{\mathcal T}(v^*_{\rm sh}\dx \ovv)-\eps\int_{-1}^\zeta\dx
\int_{z}^{\eps\zeta} v_{\rm sh}^*\dx^2 ((1+z')\ovv)\\
&-\eps\int_{-1}^\zeta\dx \int_{z}^{\eps\zeta}\dx ((1+z')\ovv)
 \dx v_{\rm sh}^*.
\end{align*}
Plugging this identity into the above expression for $A_{(3/2)}$, we
get
\begin{align}
\nonumber
A_{(3/2)}=&-2h\widetilde{\mathcal T}(v^*_{\rm sh}\dx\ovv)
+2\eps\int_{-1}^{\eps\zeta}\dx\int_{z}^{\eps\zeta}(\dx
\ovv)\dx\int_{-1}^{z'}v^*_{\rm
sh}\\
\nonumber
& -2\eps\int_{-1}^{\eps\zeta}\dx \int_z^{\eps\zeta} v^*_{\rm sh}
\dx^2 ((1+z')\ovv).
\end{align}
Using (\ref{idvs1}) together with the identities
\begin{align*}
-8h\widetilde{\mathcal T}\big(v^*_{\rm sh}\dx\ovv\big)&=-\frac{1}{3}\dx\big(h^3\dx(v^\sharp\dx\ovv)\big)-\dx\big(h^2\dx h v^\sharp \dx\ovv\big)+O(\beta),\\
-4\dx\int_{-1}^{\eps\zeta}\int_{z}^{\eps\zeta}\int_{-1}^{z'}\dx\ovv
\dx v^*_{\rm sh}&=\dx\Big(\big(\frac{h^2\dx h}{2}v^\sharp+\frac{h^3}{6}\dx v^\sharp\big)\dx \ovv\Big)+O(\beta),
\end{align*}
we finally get
\begin{equation}
A_{(3/2)}=-\frac{1}{4}\dx\big(h^3v^\sharp \dx^2\ovv\big)-\frac{1}{6}\dx\big(\dx(h^3v^\sharp)\dx\ovv\big).
\end{equation}
\subsection{Wave-current interaction in the velocity equation}\label{sect4.3}

Owing to (\ref{contribR1d}) and
(\ref{presscontr})-(\ref{presscontr2}), we deduce from the momentum
equation in (\ref{EuleravND}) that
\begin{align*}
\dsp (1+\mu{\mathcal T})\big(\dt
\ovv +\eps\ovv\dx\ovv\big)&+\dx\zeta+\eps\mu{\mathcal Q}_1(\ovv)\\
&+\eps\mu \frac{1}{h}\dx
\int_{-1+\beta b}^{\eps\zeta} \abs{v^*_{\rm sh}}^2+\mu^{3/2}\frac{1}{h}B_{(3/2)}=O(\mu^2),
\end{align*}
with
$$
B_{(3/2)}=
h\widetilde{\mathcal T}\big(\dt v^*_{\rm sh}+\eps \ovv\dx v^*_{\rm sh}\big)\\
-2\eps\int_{-1}^{\eps\zeta}\!\!\dx\!\! \int_z^{\eps\zeta} \!\! v^*_{\rm sh}
\dx^2 ((1+z')\ovv)- \dx\big(h^2 \dx\overline{v}\dx h {v^*_{\rm
  sh}}_{\vert_{z=\eps\zeta}}\big).
$$

To sum the above computations up, we have obtained the following
evolution equations for $\ovv$ (dropping
$O(\mu^2)$ terms),
\begin{equation}\label{GNvit1d}
(1+\mu{\mathcal T})\big(\dt
\ovv+\eps\ovv\dx\ovv\big)+\dx\zeta +\eps\mu {\mathcal Q}_1(\ovv)
+\eps\mu\frac{1}{h}\dx E+\eps\mu^{3/2}{\mathcal C}(v^\sharp,\ovv)=0,
\end{equation}
with
$$
E=\int_{-1+\beta b}^{\eps\zeta}\abs{v^*_{\rm
    sh}}^2
\quad\mbox{ and }\quad
v^\sharp=-\frac{24}{h^3}\int_{-1+\beta
  b}^{\eps\zeta}\int_z^{\eps\zeta}\int_{-1+\beta b}^z
v^*_{\rm sh}
$$
and ${\mathcal C}(v^\sharp,\ovv)$ as defined in (\ref{defCCb}).

When the vorticity is constant, one readily computes from
(\ref{exprVstar}) that $E= \frac{h^3}{12}\omega^2$ and
$v^\sharp=\omega h$, so that (\ref{GNvit1d}) coincides as
expected\footnote{provided that we drop in (\ref{GNvortcst}) the term
$\eps\beta\mu^{3/2}{\mathcal C}_b(\omega h,\ovv)$, which is of size
$O(\mu^2)$ under the medium amplitude bottom assumption
(\ref{assbottom}) we assumed to derive (\ref{GNvit1d}).}
with (\ref{GNvortcst}). In particular (\ref{GNvortcst}) forms a closed
system of equations in $\zeta$ and $\ovV$. This
is no longer the case in the general (variable vorticity) case since
(\ref{GNvit1d}) involves the quantity $E$ and $v^\sharp$ that need to
be determined using the vorticity equation. The presence of these
terms traduces stronger wave-current interactions than in the case of
constant vorticity.\\
The necessary closure equations on $v^\sharp$ and $E$ are derived in
the following two subsections.

\subsection{Closure equation for $v^\sharp$}\label{sect4.4}
\label{secteqvsharp1}

{\it For later use, we perform directly the
  computations in the two dimensional case $d=2$ here, the adaptations
to the case $d=1$ being straightforward.} \\
Since $v^\sharp$ appears only in the $O(\mu^{3/2})$ terms in
(\ref{GNvit1d}), it is enough to derive an equation for $v^\sharp$ at
precision $O(\sqrt{\mu})$ so that the overall $O(\mu^2)$ precision of
(\ref{GNvit1d}) is respected. \\
Such an equation is obtained by applying
the triple integration operator $
\int_{-1+\beta b}^{\eps\zeta}\!\int_{z}^{\eps\zeta}\!\!\!\int_{-1+\beta b}^{z}$
to (\ref{eqVsh}). The resulting equation is
$$
\dt (h^3 V^\sharp)+\eps \ovV\cdot \nabla (h^3V^\sharp)+\eps
h^3V^\sharp\cdot \nabla \ovV=-3 \eps h^3V^\sharp\nabla\cdot \ovV+O(\eps\sqrt{\mu}).
$$
Making use of the identity $\dt h+\eps \nabla\cdot (h\ovV)=0$, this
yields
\begin{equation}\label{eqVsharp}
\dt V^\sharp+\eps \ovV\cdot \nabla V^\sharp+\eps
V^\sharp\cdot \nabla \ovV=O(\eps\sqrt{\mu}),
\end{equation}
which is the desired closure equation.

\subsection{Closure equations for $E$}\label{sect4.5}

Since $E$ appears in an order $O(\mu)$ term in  (\ref{GNvit1d}), we need to derive an equation for $E$ at
precision $O({\mu})$ to preserve the overall $O(\mu^2)$ precision of
(\ref{GNvit1d}). The first step is to derive a more precise equation
for the shear velocity $V^*_{\rm sh}$ that takes into account the
$O(\sqrt{\mu})$ neglected in (\ref{eqVsh}). An equation for $E$ is
then deduced from this equation.

\subsubsection{An equation for $V^*_{\rm sh}$ at order $O({\mu})$}
\label{secteqvsharp2}

{\it For later use, we perform directly the
  computations in the two dimensional case $d=2$ here, the adaptations
to the case $d=1$ being straightforward.} \\
A more precise version of
(\ref{ord1}) is obtained by making more explicit the
$O(\eps\sqrt{\mu})$ term in the right-hand-side; more precisely we
substitute $O(\eps\sqrt{\mu})$ in
(\ref{ord1}) by
\begin{align*}
\eps\sqrt{\mu}\big\lbrace
&-V^*_{\rm sh}\cdot \nabla \bom_{{\mu,h}}+\big(\nabla\cdot
\int_{-1+\beta b}^z V^*_{\rm sh}\big)\dz
\bom_{\mu,h}+\bom_{\mu,h}\cdot \nabla V^*_{\rm sh}+\nabla^\perp\cdot
\ovV\dz T^*\ovV\\
&+\nabla^\perp\cdot V^*_{\rm sh}\dz V^*_{\rm sh}\big\rbrace
+O({\eps\mu}).
\end{align*}
Consequently, we include the $O(\sqrt{\mu})$ terms in (\ref{ord2}) to obtain
\begin{eqnarray*}
\dt V_{\rm sh}+\eps\ovV\cdot \nabla V_{\rm sh}+\eps V_{\rm sh}\cdot
\nabla \ovV-\eps \nabla\cdot [(1+z-\beta b)\ovV]
\dz V_{\rm sh} = \eps\sqrt{\mu} C+O(\eps\mu),
\end{eqnarray*}
with
\begin{align*}
C=&
-\int_z^{\eps\zeta} V^*_{\rm sh}\cdot \nabla \bom_{{\mu,h}}^\perp-\int_z^{\eps\zeta}(\nabla\cdot
V^*_{\rm sh})
\bom_{\mu,h}^\perp
+\Big[\big(\nabla\cdot
\int_{-1+\beta b}^z V^*_{\rm
  sh}\big)\bom_{\mu,h}^\perp\Big]_z^{\eps\zeta}\\
&
+\int_z^{\eps\zeta}\bom_{\mu,h}\cdot \nabla (V^*_{\rm
  sh})^\perp+\int_z^{\eps\zeta}\nabla^\perp\cdot \ovV(\dz
T^*\ovV)^\perp
+\int_z^{\eps\zeta}\nabla^\perp\cdot V^*_{\rm sh} \dz {V^*_{\rm
    sh}}^\perp.
\end{align*}
Using (\ref{vectorialid}) and since $\bom_{\mu,h}^\perp=-\dz V^*_{\rm
  sh}$, this yields
\begin{align*}
C=&
-\int_z^{\eps\zeta} V^*_{\rm sh}\cdot \nabla \bom_{{\mu,h}}^\perp-\int_z^{\eps\zeta}
\bom_{\mu,h}^\perp\cdot \nabla V^*_{\rm sh}
\\
&+\Big[\big(\nabla\cdot
\int_{-1+\beta b}^z V^*_{\rm
  sh}\big)\bom_{\mu,h}^\perp\Big]_z^{\eps\zeta}
+\int_z^{\eps\zeta}\nabla^\perp\cdot \ovV(\dz T^*\ovV)^\perp.
\end{align*}
Using again that $\bom_{\mu,h}^\perp=-\dz V^*_{\rm sh}$, we deduce that
\begin{align*}
C=&
\int_z^{\eps\zeta} \dz\big(V^*_{\rm sh}\cdot \nabla V^*_{\rm sh}\big)
+\Big[\big(\nabla\cdot
\int_{-1+\beta b}^z V^*_{\rm
  sh}\big)\bom_{\mu,h}^\perp\Big]_z^{\eps\zeta}
+\int_z^{\eps\zeta}\nabla^\perp\cdot \ovV(\dz T^*\ovV)^\perp\\
=&
-V^*_{\rm sh}\cdot \nabla V^*_{\rm sh}
+\ovV_{\rm sh}\cdot \nabla \ovV_{\rm sh}
+\big(\nabla\cdot
\int_{-1+\beta b}^z V^*_{\rm
  sh}\big)\dz V^*_{\rm sh}
-\nabla^\perp\cdot \ovV (T\ovV)^\perp.
\end{align*}
Finally, we have therefore the following higher order version of
(\ref{ord2}),
\begin{align*}\nonumber
\dt V_{\rm sh}+&\eps\ovV\cdot \nabla V_{\rm sh}+\eps V_{\rm sh}\cdot
\nabla \ovV +\eps\sqrt{\mu} \Big(V^*_{\rm sh}\cdot \nabla V_{\rm sh}-V_{\rm sh}\cdot
\nabla\ovV_{\rm sh}\Big)\\
\nonumber
=&\eps \nabla\cdot [\int_{-1+\beta
  b}^z(\ovV+\sqrt{\mu}V^*_{\rm sh})]
\dz V_{\rm sh}^* -\eps\sqrt{\mu}\nabla^\perp\cdot \ovV (T\ovV)^\perp+O({\eps\mu}).
\end{align*}
Integrating this yields the following equation for the shear velocity
$\ovV_{\rm sh}$,
\begin{align*}
\dt \ovV_{\rm sh}&+\eps \ovV\cdot \nabla \ovV_{\rm sh}+\eps \ovV_{\rm
  sh}\cdot \nabla \ovV+\eps\sqrt{\mu}\frac{1}{h}\nabla\cdot \int_{-1+\beta
  b}^{\eps\zeta} V^*_{\rm sh}\otimes V^*_{\rm
  sh}-\eps\sqrt{\mu}\ovV_{\rm sh}\cdot \nabla \ovV_{\rm sh}\\
&= -\eps\sqrt{\mu}\nabla^\perp\cdot \ovV (\overline{T}\ovV)^\perp+O({\eps\mu}).
\end{align*}
Taking the difference of these two equations, and dropping the
$O(\mu)$ terms, we obtain the following
higher order
generalization of (\ref{eqVsh})
\begin{eqnarray}\nonumber
\dt V^*_{\rm sh}+\eps\ovV \cdot \nabla V^*_{\rm
  sh}+\eps V_{\rm sh}^*\cdot\nabla \ovV+\eps\sqrt{\mu}\Big(V^*_{\rm sh}\cdot \nabla V^*_{\rm
  sh}-\frac{1}{h}\nabla\cdot\int_{-1+\beta b}^{\eps\zeta}V^*_{\rm sh}\otimes
V^*_{\rm sh}\Big)\\
\label{eqVshbis}
=\eps \nabla\cdot [\int_{-1+\beta
  b}^z(\ovV+\sqrt{\mu}V^*_{\rm sh})]
\dz V_{\rm sh}^*
-\eps\sqrt{\mu}\nabla^\perp\cdot \ovV (T^*\ovV)^\perp.
\end{eqnarray}

\subsubsection{An equation for $E$}

In dimension $d=1$, (\ref{eqVshbis}) takes the form
\begin{align*}\nonumber
\dt v^*_{\rm sh}+\eps\ovv \dx v^*_{\rm
  sh}+\eps v_{\rm sh}^*\dx\ovv+&\eps\sqrt{\mu}\Big(v^*_{\rm sh}\dx v^*_{\rm
  sh}-\frac{1}{h}\dx\int_{-1+\beta b}^{\eps\zeta}\abs{v^*_{\rm sh}}^2\Big)\\
%\label{eqVshbis1d}
=&\eps \dx [\int_{-1+\beta
  b}^z(\ovv+\sqrt{\mu}v^*_{\rm sh})]
\dz v_{\rm sh}^*.
\end{align*}
An equation for $E$ is simply obtained by multiplying this equation by
$v^*_{\rm sh}$ and integrating in $z$,
\begin{equation}\label{eqE}
\dt E+\eps \ovv\dx E+3\eps E \dx \ovv+\eps\sqrt{\mu}\dx F=0
\end{equation}
(up to $O(\eps\mu)$ terms), with
$$
F=\int_{-1+\beta b}^{\eps \zeta}(v^*_{\rm sh})^3.
$$
Clearly, $F$ cannot be determined in terms of $\zeta$, $\ovv$,
$v^\sharp$ and $E$, and a last equation is therefore needed.

\subsubsection{An equation for $F$}

Since $F$ appears only in the $O(\sqrt{\mu})$ term in (\ref{eqE}), we
just need to determine an evolution equation for $F$ up to
$O(\sqrt{\mu})$ terms. This equation is easily obtained by multiplying
(\ref{eqVsh}) by $\abs{v^*_{\rm sh}}^2$ and integrating vertically,
\begin{equation}
\label{eqF}
\dt F+\eps \ovv\dx F+4\eps F \dx \ovv=0
\end{equation}
(up to $O(\eps\sqrt{\mu})$ terms).

\subsection{The Green-Naghdi model}\label{sect4.6}

We are now able to give the Green-Naghdi equations in dimension $d=1$,
with general vorticity, and for non flat bottoms of medium amplitude
(i.e. $\beta=O(\sqrt{\mu})$). These equations are an order $O(\mu^2)$
approximation of (\ref{EuleravND}), where the momentum equation is
approximated by (\ref{GNvit1d}) which involves a rotational energy
$E$ determined through the finite cascade (\ref{eqE}),
(\ref{eqF}). More precisely we have, dropping $O(\mu^2)$ terms,
\begin{equation}\label{GNvort1dfin}
\begin{cases}
\dsp \dt \zeta+\dx (h\ovv)=0,\\
\dsp (1\!+\!\mu{\mathcal T})\big(\dt
\ovv\!+\!\eps\ovv\dx\ovv\big)\!+\dx\zeta +\eps\mu {\mathcal Q}_1(\ovv)
+\eps\mu\frac{1}{h}\dx E+\eps\mu^{3/2}{\mathcal C}(v^\sharp,\ovv)\!=\!0,
\\
\dsp \dt v^\sharp+\eps \ovv\dx v^\sharp+\eps
v^\sharp\dx \ovv=0,\\
\dsp \dt E+\eps \ovv\dx E+3\eps E \dx \ovv+\eps\sqrt{\mu}\dx F=0,\\
\dsp \dt F+\eps \ovv\dx F+4\eps F \dx \ovv=0.
\end{cases}
\end{equation}
where ${\mathcal T}$, ${\mathcal Q}_1$ and ${\mathcal C}$ are defined
in (\ref{defTintro}), (\ref{defQ1intro}) and (\ref{defCCb}) respectively, while we
recall that $E$, $F$ and $v^\sharp$ stand for
$$
E=\int_{-1+\beta b}^{\eps\zeta}(v^*_{\rm
    sh})^2,
\quad
F=\int_{-1+\beta b}^{\eps\zeta}(v^*_{\rm
    sh})^3
\quad
\mbox{ and }\quad
v^\sharp=-\frac{24}{h^3}\int_{-1+\beta
  b}^{\eps\zeta}\int_z^{\eps\zeta}\int_{-1+\beta b}^z
v^*_{\rm sh}.
$$
\begin{remark}\label{remNRJgen}
As in Remark \ref{remNRJcst} for the case of constant vorticity, a
local equation for the conservation of energy can be derived, which
generalizes (\ref{NRJintro}) and (\ref{NRJconstvort}), namely,
\begin{equation}\label{NRJ1d}
\dt \big({\mathfrak e}+{\mathfrak e}_{\rm rot}\big)+\dx
\big({\mathfrak F}+{\mathfrak F}_{\rm rot}\big)=0,
\end{equation}
with  (in dimensional form)
$$
{\mathfrak e}_{\rm rot}=\frac{1}{2}E\quad\mbox{ and }\quad
{\mathfrak F}_{\rm
  rot}=\frac{3}{2}E\ovv+\frac{1}{2}F+{\mathfrak F}_{{\mathcal C}}
$$
the flux ${\mathfrak F}_{{\mathcal C}}$ being as in Remark \ref{remfluxC}.
\end{remark}
\begin{remark}\label{remGNCH}
No smallness assumption on $\eps$ has been made to
derive \label{GNvort1dfin}. Assuming that $\eps=\sqrt{\mu}$ (medium
amplitude waves with the terminology of \cite{Lannes_book}), one can
simplify
(\ref{GNvort1dfin}) by dropping $O(\mu^2)$ terms.This yields
\begin{equation}\label{GNvort1dfinmed}
\begin{cases}
\dsp \dt \zeta+\dx (h\ovv)=0,\\
\dsp (1+\mu{\mathcal T})\big(\dt
\ovv+\eps\ovv\dx\ovv\big)+\dx\zeta +\eps\mu {\mathcal Q}_1(\ovv)
+\eps\mu\dx E=0,
\\
\dsp  \dt E+\eps \ovv\dx E+3\eps E \dx \ovv=0.
\end{cases}
\end{equation}
\end{remark}

\section{The $2d$ Green-Naghdi equations with general vorticity}\label{sect5}

We deal here with the derivation of Green-Naghdi type equations in the
general two-dimensional case ($d=2$). One of the main new phenomena
compared to the one-dimensional case is the interaction between the
horizontal and vertical components of the vorticity.\\
As in Section \ref{sect4}, we
assume throughout this section
that (\ref{assbottom}) holds, i.e. that the bottom variations are of
medium amplitude.

\subsection{Computation of the ``rotational Reynolds tensor''
  contribution}

Proceeding as for (\ref{contribRcst}) and (\ref{contribR1d}), we
write, up to $O(\mu^2)$ terms,
\begin{align*}
\nonumber
\ep\mu \nabla\cdot \int_{-1+\beta b}^{\ep\zeta}V^*\otimes
V^*=&\ep\mu\nabla\cdot\int_{-1+\beta b}^{\ep\zeta}\vs\otimes\vs\\
&-\ep\mu^{3/2}\frac{1}{2}
\nabla\cdot\int_{-1+\beta b}^{\ep\zeta}\int_z^{\eps\zeta}
\nabla\nabla\cdot \big((1+z-\beta b)\ovV\big)\otimes \vs\\
&-\ep\mu^{3/2}\frac{1}{2}
\nabla\cdot\int_{-1+\beta b}^{\ep\zeta}\int_z^{\eps\zeta}\vs\otimes\nabla\nabla\cdot \big((1+z-\beta b)\ovV\big).
\end{align*}
Introducing
\begin{eqnarray*}
V^\sharp=\frac{24}{h^3}\int_{-1+\beta b}^{\ep \zeta}\int_{z}^{\ep \zeta}\int_{-1+\beta b}^z \vs.
\end{eqnarray*}
and proceeding as in \S \ref{sect4.1}, we obtain
\begin{align}
\nonumber
\ep\mu \nabla\!\cdot\! \int_{-1+\beta b}^{\ep\zeta}\!\!\!V^*\!\otimes\!
V^*=&\ep\mu\nabla\cdot\int_{-1+\beta b}^{\ep\zeta}\vs\otimes\vs\\
\label{contribR2d}
&-\eps\mu^{3/2}
\frac{1}{24}\nabla\cdot\big(h^3(V^\sharp\otimes\nabla\nabla\!\cdot\!\ovV+\nabla\nabla\!\cdot\!
\ovV\otimes V^\sharp)\big)+O(\mu^2).
\end{align}
\subsection{Computation of the pressure contribution}\label{sect5.2}

As for (\ref{contribPcst}) and (\ref{presscontr}), we have
\begin{align}
\nonumber
\frac{1}{\ep}\int_{-1+\beta b}^{\eps\zeta}\nabla p &= h\nabla\zeta
+\int_{-1+\beta b}^{\ep \zeta}\nabla\int_{z}^{\ep \zeta}\left(\pa_t w
  +\ep V\cdot \nabla w +\frac{\ep}{\mu}w\pa_z w\right)\\
\label{vac2}
&= h\nabla\zeta+ \mu A_{(1)}+\mu^{3/2} A_{(3/2)}+O(\mu^2),
\end{align}
where, as in the irrotational case,
\begin{eqnarray}
\label{vac3}
A_{(1)}= h \mathcal{T}\big(\pa_t \ovV+\eps \ovV\cdot \nabla
\ovV\big)+\ep h {\mathcal Q}_1(\ovV),
\end{eqnarray}
and ${\mathcal Q}_1(\ovV)$ as in (\ref{defQ1intro}). For the $O(\mu^{3/2})$
component, we have
\begin{align*}
\nonumber
A_{(3/2)}
=&h\widetilde{\mathcal T}(\dt V^*_{\rm
  sh}\!+\!\eps \ovV\cdot \nabla V^*_{\rm sh})
+2\eps\int_{-1}^{\eps\zeta}\nabla\int_{z}^{\eps\zeta}\big(\int_{-1}^z\nabla\cdot
\vs\big)\nabla\cdot \ovV
\\
\nonumber
&\!-\!\eps \int_{-1}^{\eps\zeta}\!\!\nabla\!\! \int_z^{\eps\zeta} \!\!V^*_{\rm sh}
\cdot \nabla \big(\nabla\cdot ( (1+z')\ovV)\big)+\eps\int_{-1}^{\eps\zeta}\nabla\int_z^{\eps\zeta}\nabla\cdot
\big((1+z)\ovV\big)\nabla\cdot \vs\\
&+\eps
\int_{-1}^{\eps\zeta}\nabla\int_{z}^{\eps\zeta}\big([\dive,\ovV\cdot
\nabla]-(\nabla\cdot \ovV)\nabla^T\big)\int_{-1}^z  V^*_{\rm sh}.
\end{align*}
Applying
$\widetilde{\mathcal T}$ to (\ref{eqVsh}) to handle the first term of
the right-hand-side, we obtain proceeding as in \S \ref{sectcontribP1d} that
\begin{align}
\nonumber
A_{(3/2)}
=&-\eps h\widetilde{\mathcal T}(V^*_{\rm
  sh}\cdot \nabla\ovV+ V^*_{\rm sh} \nabla\cdot \ovV)
+2\eps\int_{-1}^{\eps\zeta}\nabla\int_{z}^{\eps\zeta}\big(\int_{-1}^z\nabla\cdot
\vs\big)\nabla\cdot \ovV
\\
\nonumber
&\!-\!2\eps \int_{-1}^{\eps\zeta}\!\!\nabla\!\! \int_z^{\eps\zeta} \!\!V^*_{\rm sh}
\cdot \nabla \big(\nabla\cdot ( (1+z')\ovV)\big)\\
\nonumber
&+\eps
\int_{-1}^{\eps\zeta}\nabla\int_{z}^{\eps\zeta}\big([\dive,\ovV\cdot
\nabla]-(\nabla\cdot \ovV)\nabla^T\big)\int_{-1}^z V^*_{\rm sh}\\
\label{eqexam}
=&-\eps\nabla \big[\frac{h^3}{4}\vs\cdot \nabla\nabla\cdot
\ovV+\frac{1}{12}\nabla\cdot \ovV\nabla\cdot (h^3\vs)+\frac{1}{12}{\rm
  Tr}\big(\nabla \ovV\nabla(h^3V^\sharp)\big)\big].
\end{align}

\subsection{Wave-current interaction in the velocity equation}

From the momentum
equation in (\ref{EuleravND})  and (\ref{contribR2d}), (\ref{vac2}),
(\ref{vac3}) and (\ref{eqexam}), we obtain

\begin{eqnarray*}
\left(1+\mu\mathcal{T}\right)\big(\pa_t\ovV +\ep \ovV\cdot
\nabla\ovV\big)+\mu \ep  {\mathcal Q}_1(\ovV)+\eps\mu\frac{1}{h}\nabla\cdot E+\eps\mu^{3/2}{\mathcal C}(V^\sharp,\ovV)=0,
\end{eqnarray*}
where
$$
E=\int_{-1+\beta b}^{\eps\zeta}\vs\otimes \vs
\quad\mbox{ and }\quad
V^\sharp=-\frac{24}{h^3}\int_{-1+\beta
  b}^{\eps\zeta}\int_z^{\eps\zeta}\int_{-1+\beta b}^z
\vs,
$$
and where ${\mathcal C}(V^\sharp,\ovV)$ is the two-dimensional
generalization of (\ref{defCCb}),
\begin{align}
\nonumber
{\mathcal C}(V^\sharp,\ovV)=&-\frac{1}{24h}\nabla\cdot\big(h^3(V^\sharp\otimes\nabla\nabla\cdot\ovV+\nabla\nabla\cdot
\ovV\otimes V^\sharp)\big) \\
\label{defC2d}
&-\frac{1}{4h}\nabla \big[h^3V^\sharp\cdot \nabla\nabla\cdot
\ovV+\frac{1}{3}\nabla\cdot \ovV\nabla\cdot (h^3V^\sharp)+\frac{1}{3}{\rm
  Tr}\big(\nabla \ovV\nabla(h^3V^\sharp)\big)\big].
\end{align}

As in the one dimensional case (see \S \ref{sect4.3}), closure
equations are needed for $V^\sharp$ and for $E$. The closure equation
has already been derived in the two-dimensional case in \S
\ref{sect4.4}; it is given by
\begin{eqnarray*}
\pa_t V^\sharp +\ep (\ovV\cdot\nabla)V^\sharp+\ep(V^\sharp\cdot\nabla)\ovV=0.
\end{eqnarray*}
The derivation of closure equations for $E$ is addressed in the
following section.

\subsection{Closure equation for E}

We first derive an equation for $E$, which involves a third order
tensor $F$, for which we also derive an equation that closes the system.
\subsubsection{An equation for $E$}
Recall first that we obtained in (\ref{eqVshbis}) an evolution equation
for $\vs$, namely
\begin{eqnarray*}\nonumber
\dt V^*_{\rm sh}&+\eps\ovV \cdot \nabla V^*_{\rm
  sh}+\eps V_{\rm sh}^*\cdot\nabla \ovV
=\eps \nabla\cdot \big((1+z-\beta b)\ovV\big)\dz\vs+\eps\sqrt{\mu} C^*
\end{eqnarray*}
with
$$
C^*=-V^*_{\rm sh}\cdot \nabla V^*_{\rm
  sh}+\frac{1}{h}\nabla\cdot\int_{-1+\beta b}^{\eps\zeta}V^*_{\rm sh}\otimes
V^*_{\rm sh}+\nabla\cdot (\int_{-1+\beta
  b}^zV^*_{\rm sh})
\dz V_{\rm sh}^*
-\nabla^\perp\cdot \ovV (T^*\ovV)^\perp.
$$
Time differentiating the tensor $E$, we get therefore
\begin{align*}
\dt E=&\dt (\eps\zeta)
(\vs\otimes\vs)_{\vert_{\rm surf}} +\int_{-1+\beta b}^{\eps\zeta}
\left(\pa_t\vs\otimes\vs+\vs\otimes\pa_t\vs\right)\\
= &\dt (\eps\zeta)
(\vs\otimes\vs)_{\vert_{\rm surf}} +I_1+I_2+I_3+\eps\sqrt{\mu}G
\end{align*}
with
\begin{align*}
I_1=&-\eps\int_{-1+\beta
  b}^{\eps \zeta} \left((\ovV\cdot \nabla)\vs \otimes \vs+\vs\otimes
  (\ovV\cdot\nabla)\vs\right)\\
=&-\eps \ovV\cdot \nabla E+\eps\big(\ovV\cdot
\nabla(\eps\zeta)(\vs\otimes\vs)_{\vert_{\rm surf}} -\ovV\cdot
\nabla(\beta b)(\vs\otimes\vs)_{\vert_{\rm bott}} \big),
\end{align*}
while
\begin{align*}
I_2=&\int_{-1+\beta b}^{\eps\zeta} \nabla\cdot \big((1+z-\beta
b)\ovV\big)\left(\pa_z \vs \otimes\vs +\vs\otimes \pa_z \vs\right)\\
=&-E+h (\vs\otimes\vs)_{\vert_{\rm surf}} -\nabla(\beta b)\cdot \ovV\big((\vs\otimes\vs)_{\vert_{\rm surf}} -(\vs\otimes\vs)_{\vert_{\rm bott}} \big),
\end{align*}
and
\begin{align*}
I_3=&-\eps\int_{-1+\beta b}^{\eps\zeta}
\left((\vs\cdot\nabla\ovV)\otimes \vs + \vs\otimes
  (\vs\cdot\nabla)\ovV\right)\\
=&-(\nabla\ovV)^\perp E- E\nabla\ovV.
\end{align*}
Finally, the matrix $G$ is given by
\begin{align*}
G=&\int_{-1+\beta b}^{\eps \zeta} C^*\otimes \vs+\vs\otimes C^*\\
=&-\!\!\!\int_{-1+\beta b}^{\eps \zeta}\!\!\!(\vs\cdot
\nabla\vs)\otimes\vs+\vs\otimes (\vs\cdot \nabla\vs)\\
&-\!\!\!\int_{-1+\beta b}^{\eps\zeta}\!\!\! (\nabla\cdot \vs)\vs\otimes
\vs-\eps\nabla\zeta\cdot \vs
(\vs\otimes\vs)_{\vert_{\rm surf}}
+\beta \nabla  b\cdot \vs   (\vs\otimes\vs)_{\vert_{\rm bott}}\\
&+\frac{h^4}{24}\big(\nabla^\perp\cdot \ovV\nabla^\perp\nabla\cdot
\ovV\otimes V^\sharp+V^\sharp\otimes \nabla^\perp\cdot \ovV\nabla^\perp\nabla\cdot
\ovV\big).
\end{align*}
Gathering all these computations, we finally get that
\begin{eqnarray}
\pa_t E + \ep \ovV\cdot\!\nabla \! E +\ep \nabla\!\cdot\!\ovV E +\ep \nabla
\ovV ^T E+\ep E\nabla\ovV+\eps\sqrt{\mu}\nabla\!\cdot\! F
\label{eqE2d}
=\eps\sqrt{\mu}
{\mathcal D}(V^\sharp,\ovV)
\end{eqnarray}
where
\begin{eqnarray*}
F=\int_{-1+\beta b}^{\ep\zeta} \vs\otimes \vs\otimes \vs
\end{eqnarray*}
and
\begin{eqnarray*}
{\mathcal D}(V^\sharp,\ovV) =\frac{h^3}{24}\nabla^\perp\cdot \ovV \big(\nabla^\perp \nabla\cdot\ovV\otimes
V^\sharp+V^\sharp\otimes  \nabla^\perp \nabla\cdot\ovV \big).
\end{eqnarray*}
\begin{remark}
The tensor $E$ is symmetric, and therefore $E=E^T$. The equation
(\ref{eqE2d}) is therefore the same as the equation governing the
evolution of the Reynolds tensor in barotropic turbulent compressible
fluids (\cite{MP,Pope,GG}) with a source term ${\bf S}=-\eps\sqrt{\mu}\nabla \cdot  F+\eps\sqrt{\mu}
{\mathcal D}(V^\sharp,\ovV)$.  The structure of the source term ${\bf
  S}$ is central in turbulence theory and is still under intense
investigation; in the present case, it has a well defined structure
and the system of equations can be closed by
deriving evolution equations on $V^\sharp$ and $F$. We refer to \cite{GG}
for an instructive geometric study of the equation in the case ${\bf S}=0$.
\end{remark}
\subsubsection{An equation for $F$}

To close the system we need and equation for the third order tensor $F$ up to order
$O(\sqrt{\mu})$. We just compute, for $i,j,k=1,2$,
\begin{eqnarray*}
\pa_t F_{ijk}= \eps\dt\zeta {(\vs)_i(\vs)_j(\vs)_k}_{\vert_{\rm
    surf}}+\int_{-1+\beta b}^{\eps \zeta}
\pa_t\big((\vs)_i(\vs)_j(\vs)_k\big)
\end{eqnarray*}
and use (\ref{eqVsh}) to get, up to $O(\sqrt{\mu})$ terms and with
Einstein's summation convention on repeated indices,
\begin{align*}
\int_{-1+\beta b}^{\eps \zeta} \pa_t (\vs)_i (\vs)_j(\vs)_k=&-\eps \int_{-1+\beta b}^{\eps \zeta}  \ovV\cdot\nabla
((\vs)_i(\vs)_j(\vs)_k)\\
&-F_{lkj}\pa_l\ovV_i-F_{ilk}\pa_l\ovV_j-F_{ijl}\pa_l\ovV_k\\
&+\eps\int_{-1+\beta b}^{\eps\zeta} \nabla\cdot\big((1+z-\beta b)\ovV\big)\pa_z((\vs)_i(\vs)_j(\vs)_k).
\end{align*}
One then readily deduces the following equation
\begin{eqnarray}
\label{eqF2d}
\pa_t F_{ijk}+\ep \ovV\cdot\nabla F_{ijk}+\ep F_{lkj}\pa_l \ovV_i+\ep F_{ilk}\pa_l\ovV_j+\ep F_{ijl}\pa_l\ovV_k+\ep \nabla\cdot \ovV F_{ijk}=0
\end{eqnarray}
($i,j,k=1,2$), which closes the system of equations on $\zeta$, $\ovV$, $V^\sharp$,
$E$ and $F$.

\subsection{The Green-Naghdi model}\label{sectGNmodel2D}

We can now give the two-dimensional generalization of the Green-Naghdi
equations (\ref{GNvort1dfin}) with general vorticity. More precisely we have, dropping $O(\mu^2)$ terms,
\begin{equation}\label{GNvort2dfin}
\begin{cases}
\dsp\dt \zeta+\nabla\cdot (h\ovV)=0,\\
\dsp \left(1\!+\!\mu\mathcal{T}\right)\big(\pa_t\ovV \!+\!\ep \ovV\!\cdot\!
\nabla\ovV\big)+\nabla\zeta+\mu \ep  {\mathcal Q}_1(\ovV)+\eps\mu\frac{1}{h}\nabla\!\cdot\! E+\eps\mu^{3/2}{\mathcal C}(V^\sharp,\ovV)=0,\\
\dsp \pa_t V^\sharp +\ep
(\ovV\cdot\nabla)V^\sharp+\ep(V^\sharp\cdot\nabla)\ovV=0,\\
\dsp \pa_t E + \ep \ovV\cdot\!\nabla \! E +\ep \nabla\!\cdot\!\ovV E +\ep \nabla
\ovV ^T E+\ep E\nabla\ovV+\eps\sqrt{\mu}\nabla\!\cdot\! F
=\eps\sqrt{\mu}
{\mathcal D}(V^\sharp,\ovV) \\
\dsp\pa_t F_{ijk}+\ep \ovV\cdot\nabla F_{ijk}+\ep F_{lkj}\pa_l
\ovV_i+\ep F_{ilk}\pa_l\ovV_j+\ep F_{ijl}\pa_l\ovV_k+\ep \nabla\cdot
\ovV F_{ijk}=0,
\end{cases}
\end{equation}
($i,j,k=1,2$),
where ${\mathcal T}$, ${\mathcal Q}_1$ and ${\mathcal C}$ are defined
in (\ref{defTintro}), (\ref{defQ1intro}) and (\ref{defC2d}) respectively, while we
recall that $E$, $F$ and $V^\sharp$ stand for
$$
E=\int_{-1+\beta b}^{\eps\zeta}\vs\otimes \vs,
\qquad
F=\int_{-1+\beta b}^{\eps\zeta}\vs\otimes\vs\otimes\vs
$$
and
$$
V^\sharp=-\frac{24}{h^3}\int_{-1+\beta
  b}^{\eps\zeta}\int_z^{\eps\zeta}\int_{-1+\beta b}^z
V^*_{\rm sh}.
$$
and where
\begin{equation}\label{defD}
{\mathcal D}(V^\sharp,\ovV) =\frac{h^3}{24}\nabla^\perp\cdot \ovV \big(\nabla^\perp \nabla\cdot\ovV\otimes
V^\sharp+V^\sharp\otimes  \nabla^\perp \nabla\cdot\ovV \big).
\end{equation}
\begin{remark}\label{remfluxC2d}
In dimension $d=2$, the structural property given in Remark
\ref{remfluxC} can be generalized; there holds
$$
h{\mathcal C}(V^\sharp,\ovV)\cdot \ovV=\nabla\cdot {\mathfrak
  F}_{\mathcal C}+\frac{1}{2}\mbox{Tr}({\mathcal D}(V^\sharp,\ovV)),
$$
where the flux ${\mathfrak F}_{\mathcal C}$ is given by the vector
\begin{align*}
 {\mathfrak F}_{\mathcal C}=&-\frac{h^3}{12}\Big( \frac{1}{2}V^\sharp
 \cdot \ovV(\nabla\nabla\cdot V)+\frac{1}{2}\ovV\cdot
 \nabla\nabla\cdot V V^\sharp+3V^\sharp\cdot \nabla \nabla\cdot \ovV
 \ovV-(\nabla\cdot \ovV)^2V^\sharp\Big)\\
&+\frac{1}{12}\Big(\nabla\cdot \ovV\nabla\cdot
(h^3V^\sharp)\ovV+\mbox{Tr}\big(\nabla\ovV\nabla(h^3V^\sharp)\big)\ovV\Big)\\
&+\frac{h^3}{12}\Big((\nabla\cdot \ovV)(V_1^\sharp\nabla^\perp
\ovV_2-V_2^\sharp\nabla^\perp \ovV_1)\Big).
\end{align*}
Contrary to the one dimensional case, the quantity $h{\mathcal
  C}(V^\sharp,\ovV)\cdot \ovV$ is not the divergence of a flux. The
presence of the term $\frac{1}{2}\mbox{Tr}({\mathcal
  D}(V^\sharp,\ovV))$ traduces a new mechanism of energy flux due to
the interaction of the vertical and horizontal components of the vorticity.
\end{remark}
\begin{remark}\label{remNRJ2d}
Using Remark \ref{remfluxC2d}, the local equation for the conservation of energy derived in Remark
\ref{remNRJgen} in one dimension can be generalized in two dimensions;
one has
\begin{equation}\label{NRJ2d}
\dt \big({\mathfrak e}+{\mathfrak e}_{\rm rot}\big)+\nabla\cdot
\big({\mathfrak F}+{\mathfrak F}_{\rm rot}\big)=0,
\end{equation}
with  (in dimensional form)
$$
{\mathfrak e}_{\rm rot}=\frac{1}{2}\mbox{Tr}\,E\quad\mbox{ and }\quad
{\mathfrak F}_{\rm
  rot}=\frac{1}{2}(\mbox{Tr}\,E)\ovV+E\ovV+\frac{1}{2}\left(\begin{array}{c}F^{111}+F^{122}\\
F^{211}+F^{222}\end{array}\right)+{\mathfrak F}_{{\mathcal C}}
$$
the flux ${\mathfrak F}_{{\mathcal C}}$ being as in Remark \ref{remfluxC2d}.
\end{remark}

\section{Reconstruction of the velocity profile and vorticity dynamics}\label{sect6}

\subsection{First order (St-Venant) reconstruction of the velocity}\label{sect6.1}

As noted in Remark \ref{remSV} a first order approximation of the
averaged Euler equations (\ref{EuleravND}) is provided by the
Saint-Venant (or nonlinear shallow water) equations
$$
\begin{cases}
\dsp \dt \zeta+\nabla\cdot (h\ovV)=0,\\
\dsp \dt (h\ovV)+h\nabla\zeta+\eps \nabla\cdot (h\ovV\otimes \ovV)=0;
\end{cases}
$$
these equations are the same as in the irrotational case. However,
when reconstructing the velocity fluid ${\bf V}(t,X,z)$ inside the
fluid domain, the effects of the vorticity cannot be
neglected. According to (\ref{Vmu}) and (\ref{wmu}), the horizontal and vertical velocities
are given at first order by
$$
{\bf V}=\ovV+\sqrt{\mu}\vs,\quad\mbox{ and }\quad
\frac{1}{\mu}{\bf w}=-\nabla\cdot \big[(1+z-\beta b)\ovV
\big]-\mu^{1/2}\nabla\cdot \int_{-1+\beta b}^{z}V^*_{\rm sh}.
$$
We also know that $\vs$ is found through the resolution of
(\ref{eqVsh}). This equation is a variable coefficients linear
equation cast on the fluid domain $\Omega_t$. Though $\Omega_t$ is
known at this step (through the resolution of the St-Venant equation),
it is still a moving, $d+1$ dimensional domain and the numerical
computation of the solutions to (\ref{eqVsh}) is time consuming. We
propose here a simpler approach consisting in deriving a simple,
scalar, $d$ dimensional equation determining the horizontal velocity on
each level line $\Gamma_\theta$, with
$$
\Gamma_{\theta}=\{ (X,z), z=-1+\beta
b(X)+\theta h(t,X)\}\qquad (\theta\in [0,1]),
$$
so that $\Gamma_0$ corresponds to the bottom and $\Gamma_1$ to the
surface. With the notation
\begin{equation}\label{defnotaVtheta}
V^*_{\theta}(X)={\vs}_{\vert_{\Gamma_\theta}}(X)=\vs(X,-1+\beta
b(X)+\theta h(t,X)),
\end{equation}
one has
\begin{equation}\label{recVw}
\begin{array}{lcl}
\dsp {\bf V}_{\vert_{\Gamma_\theta}}&=&\dsp \ovV+\sqrt{\mu}V^*_\theta,\\
\dsp \frac{1}{\mu}{\bf w}_{\vert_{\Gamma_\theta}}&=&\dsp -\nabla\cdot \big(h
(\theta\ovV+\sqrt{\mu}Q_\theta)\big)+\nabla(-1+\beta b+\theta h)\cdot
(\ovV+\sqrt{\mu}V^*_\theta),
\end{array}
\end{equation}
with
$$
Q_\theta=\frac{1}{h}\int_{-1+\beta b}^{-1+\beta b+\theta
  h}\vs dz=\int_0^\theta V^*_{\theta'} d\theta';
$$
and one easily gets from (\ref{eqVsh}) that
\begin{equation}\label{eqVtheta}
\dt V^*_\theta+\eps\ovV\cdot \nabla V_\theta^*+\eps V_\theta^*\cdot \nabla \ovV=0.
\end{equation}
and that
\begin{equation}\label{eqQtheta}
\dt Q_\theta+\eps\ovV\cdot \nabla Q_\theta+\eps Q_\theta\cdot \nabla \ovV=0.
\end{equation}
The velocity field in the fluid domain can therefore be fully determined at
the precision of the model by the following {\it decoupled} procedure:
\begin{enumerate}
\item Solve the Saint-Venant equations (\ref{StV}) to get $\zeta$ and
  $\ovV$ on the desired time interval.
\item Get the quantities $V_\theta^*$ and $Q_\theta$ from their
  initial values by solving (\ref{eqVtheta})-(\ref{eqQtheta}).
\item Reconstruct the velocity field on each level line
  $\Gamma_\theta$ ($0\leq \theta\leq 1$) by
  using the formulas (\ref{recVw}).
\end{enumerate}
Though the first order (St-Venant) approximation is the same as in the
irrotational case, the velocity in the fluid domain therefore differs
from the irrotational theory by the shear component
$\sqrt{\mu}V^*_\theta$ which is essentially advected at the mean
velocity $\ovV$.

\subsection{Second order (Green-Naghdi) reconstruction of the velocity}\label{sect6.2}

As seen in \S \ref{sectGNmodel2D}, a second order approximation of the
averaged Euler equations (\ref{EuleravND}) is provided by the
Green-Naghdi equations with vorticity (\ref{GNvort2dfin}). According
to (\ref{Vmu2}) and (\ref{wmu2}), the velocity in the fluid domain is
given by
\begin{equation}\label{eqVGiro}
{\bf V}=\ovV+\sqrt{\mu}V^*_{\rm sh}    +\mu T^*\ovV +O(\mu^{3/2})
\end{equation}
and 
$$
\frac{1}{\sqrt{\mu}}{\bf w}=-\nabla\cdot \big[(1+z-\beta b)\ovV
\big]-\mu^{1/2}\nabla\cdot \int_{-1+\beta b}^{z}V^*_{\rm sh}-\mu\nabla\cdot \int_{-1+\beta b}^{z}T^*\ovV+O(\mu^{3/2}).
$$
Using the same notations as in the previous section, and defining
$T^*_\theta$ as
$$
T^*_\theta\ovV=-\frac{1}{2}(\theta^2-\frac{1}{3})h^2\nabla\nabla\cdot
\ovV+\beta (\theta-\frac{1}{2})h \big( \nabla b\cdot \nabla
\ovV+\nabla(\nabla b\cdot \ovV)\big),
$$
we therefore get up to $O(\mu^{3/2})$ terms\footnote{It is actually
  possible to reconstruct the velocity at order $O(\mu^2)$, but this
  requires to include the $O(\mu)$ termes in the equations for
  $V_\theta^*$; for the sake of simplicity, we therefore stick to
  order $O(\mu^{3/2})$.},
\begin{equation}\label{recVw2}
\begin{array}{lcl}
\dsp {\bf
  V}_{\vert_{\Gamma_\theta}}&=&\dsp \ovV+\sqrt{\mu}V^*_\theta+\mu T^*_\theta
\ovV,\\
\dsp \frac{1}{\mu}{\bf w}_{\vert_{\Gamma_\theta}}&=&\dsp -\nabla\cdot \big(h
(\theta\ovV+\sqrt{\mu}Q_\theta)\big)+\nabla(-1+\beta b+\theta h)\cdot
(\ovV+\sqrt{\mu}V^*_\theta)\\
& &-\mu\nabla\big(h\cdot \int_{0}^\theta
T^*_{\theta'}\ovV d\theta'\big)+\mu\nabla(-1+\beta b+\theta h)T_\theta^*\ovV.
\end{array}
\end{equation}
We also deduce from (\ref{eqVshbis}) that
\begin{equation}
\dt V^*_\theta+\ovV\cdot \nabla V_\theta^*+V_\theta^*\cdot \nabla
\ovV+\eps\sqrt{\mu}V_\theta^*\cdot \nabla
V_\theta^*=\eps\sqrt{\mu}({\mathcal S}_1+{\mathcal S}_2),
\label{eqVthetabis}
\end{equation}
where the source terms ${\mathcal S}_1$ and ${\mathcal S}_2$ are given by
\begin{align*}
{\mathcal S}_1&=\frac{1}{h}\nabla\cdot E+\frac{q_\theta}{h}\nabla\cdot
(hQ_\theta),\\
{\mathcal S}_2&=
-\frac{h^2}{6}(1-3\theta^2)\nabla^\perp\cdot
\ovV\nabla^\perp\nabla\cdot \ovV,
\end{align*}
the quantities $q_\theta$ and $Q_\theta$ being defined as
$$
q_\theta=h\dz {\vs}_{\vert_{\Gamma_\theta}}=\partial_\theta V_\theta^*
\quad\mbox{ and }\quad
Q_\theta=\frac{1}{h}\int_{-1+\beta b}^{-1+\beta b+\theta
  h}\vs dz=\int_0^\theta V^*_{\theta'} d\theta';
$$
these two quantities can be straightforwardly computed from their
initial values by solving the equations
\begin{equation}\label{eqQq}
\begin{array}{lcl}
\dsp \dt q_\theta +\eps q_\theta\cdot \nabla \ovV+\eps \ovV\cdot \nabla
q_\theta&=&0,\\
\dsp \dt Q_\theta +\eps Q_\theta\cdot \nabla \ovV+\eps \ovV\cdot
\nabla Q_\theta&=&0.
\end{array}
\end{equation}

The velocity field in the fluid domain can therefore be fully determined at
the precision of the model by the following {\it decoupled} procedure:
\begin{enumerate}
\item Solve the Green-Naghdi equations (\ref{GNvort2dfin}) to get $\zeta$ and
  $\ovV$ on the desired time interval.
\item Get the quantities $V_\theta^*$, $q_\theta$ and $Q_\theta$ from their
  initial values by solving (\ref{eqVthetabis}) and (\ref{eqQq}).
\item Reconstruct the velocity field on each level line
  $\Gamma_\theta$ ($0\leq \theta\leq 1$) by
  using the formulas (\ref{recVw2}).
\end{enumerate}
There are several important differences to be underlined if one compares
this second order approximation to the first-order (St-Venant)
approximation considered in the previous section:
\begin{itemize}
\item The equation for $\zeta$ and $\ovV$ are not the same as in the
  irrotational theory;
\item The quadratic term in the left-hand-side and the source term ${\mathcal S}_1$ in (\ref{eqVthetabis}) traduce
  a more complex behavior of the shear velocity due to the horizontal vorticity.
\item The above reconstruction procedure exhibits a
  mechanism of creation of shear velocity from vertical
  vorticity. Even if we start from an initial zero horizontal
  vorticity (and therefore ${\vs}_{\vert_{t=0}}=0$), the shear
  velocity does not remain equal to zero during the evolution of the
  flow. Indeed, due to the source term ${\mathcal S}_2$ in
  (\ref{eqVthetabis}) (it self proportional to the vertical vorticity
  $\nabla^\perp\cdot \ovV$), the quantity $V_\theta^*$ departs from
  its zero initial value.
\end{itemize}

\subsection{The dynamics of the vertical vorticity}\label{sect6.3}

The dynamics of the vertical component $\bom_{\mu,v}=\nabla^\perp
\cdot {\bf V}$ of the vorticity is very important for the study of
rip currents for instance \cite{Hammack,Chen1}). One directly gets from (\ref{eqVGiro}) that,
up to $O(\mu^{3/2})$ terms,
$$
\bom_{\mu,v}=\nabla^\perp\cdot \ovV+\sqrt{\mu}\nabla^\perp \cdot \vs+\mu
\frac{h}{3}\nabla^\perp h\cdot \nabla (\nabla\cdot \ovV).
$$
In order to make more explicit the possible creation of vertical
vorticity from horizontal vorticity, we consider here the time
evolution of the vertically averaged vertical velocity
$\overline{\bom}_{\mu,v}$ defined as
$$
\overline{\bom}_{\mu,v}(t,X)=\frac{1}{h}\int_{-1+\beta b}^{\eps\zeta}\bom_{\mu,v}(t,X,z)dz
$$
so that
\begin{equation}\label{avvort}
\overline{\bom}_{\mu,v}(t,X)=\overline{\omega}_0+\sqrt{\mu}\overline{\omega}_1,
\end{equation}
with
\begin{equation}\label{eqvortv1}
\begin{array}{lcl}
\dsp \overline{\omega}_0&=&\dsp \nabla^\perp\cdot \ovV+\mu
\frac{h}{3}\nabla^\perp h\cdot \nabla (\nabla\cdot \ovV)\\
\dsp \overline{\omega}_1&=&\dsp -\frac{1}{h}\big(\eps\nabla^\perp \zeta\cdot V^*_1-\beta
\nabla^\perp bV^*_0\big)
\end{array}
\end{equation}
(according to the notation (\ref{defnotaVtheta}), $V_0^*$ and $V_1^*$
correspond to the evaluation of $\vs$ at the bottom and at the surface
respectively). We have already seen that all the $V^*_\theta$ can be
recovered from their initial value through (\ref{eqVthetabis});
considering the particular cases $\theta=0,1$ we therefore get
$\overline{\omega}_1$. For $\overline{\omega}_0$, we apply
$\nabla^\perp$ to the second equation of (\ref{GNvort2dfin}) so that
\begin{equation}\label{eqvortv0}
\dt \overline{\omega}_0+\nabla\cdot
(\overline{\omega}_0\overline{V})+\eps\mu \nabla^\perp\cdot
\big(\frac{1}{h}\nabla\cdot E\big)+\eps\mu^{3/2}\nabla^\perp\cdot
{\mathcal C}(V^\sharp,\ovV)=0.
\end{equation}
The procedure to recover the vertical vorticity is therefore the
following:
\begin{enumerate}
\item Solve the Green-Naghdi equations (\ref{GNvort2dfin}) to get
  $\zeta$, $\ovV$, $E$ and $V^\sharp$ on the desired time interval.
\item Get $\overline{\omega}_0$ from its initial value by solving (\ref{eqvortv0}).
\item Get the quantities $V_0^*$ and $V_1^*$ from their
  initial values by solving (\ref{eqVthetabis}) and deduce
  $\underline{\omega}_1$ from (\ref{eqvortv1}).
\item Recover the averaged vorticity $\overline{\bom}_{\mu,v}$ through (\ref{avvort}).
\end{enumerate}
One important aspect to underline here is that it is possible to start
with a zero averaged vertical vorticity but that this quantity becomes
nonzero with the evolution of the flow provided that some (horizontal)
vorticity is present. This mechanism of transfer from the horizontal
to the vertical vorticity is likely to play an important role in the
study of rip-currents; its detailed study is left for future works.
\section{Conclusion}\label{sect7}

We have derived here several fully nonlinear models generalizing the
Green-Naghdi equations in presence of vorticity; we have based our
formal computations on the rigorous estimates derived in
\cite{CastroLannes}, hereby ensuring the validity of the
approximations made throughout this article. For the sake of
clarity, we have presented these models by increasing complexity
(constant vorticity in dimension $d=1$, general vorticity in dimension
$d=1$, general vorticity in dimension $d=2$). The most remarkable
feature of these models is that they do not require the coupling with
a $d+1$-dimensional equation for the vorticity; indeed, they differ from
the irrotational theory by the coupling, reminiscent of turbulence theory, with a {\it finite} cascade of
equations. This cascade gives the evolution of the components of the ``Reynolds'' tensor
describing the self-interaction of the shear velocity induced by the
vorticity on the one hand, and its interaction with the
non-hydrostatic vertical variations typical of Boussinesq type models
on the other hand. The reconstruction of
the velocity profile in the ($d+1$ dimensional) fluid domain can then
be performed by a completely decoupled and simple procedure. The most
striking phenomena here are a mechanism of {\it creation of horizontal
  shear} from vertical vorticity and conversely, a mechanism of
{\it transfer of horizontal to vertical vorticity} likely to play a
key-role in rip-currents for instance.\\
 A natural
perspective opened by this work is therefore to allow for the
possibility of shear flows in the numerous codes developed recently
for the numerical simulation of (irrotational) Green-Naghdi systems
(see for instance
\cite{CKDKC,Cienfuegos,LGH,BCLMT,Delis,Dutykh,Mario,MarcheLannes} and
the review \cite{Betal}) and the modeling of rip-currents using the
models derived in this article.

In view of describing rip-currents, another natural
perspective is to take into account the creation of vorticity by
mechanisms such as wave breaking, surface and bottom boundary layers,
etc. The present work deals indeed with a conservative framework (as
shown for instance by the equation (\ref{NRJ2d}) for the local
conservation of energy). Our goal here was to understand the coupling
between surface waves and underlying vortical flows. It complements in
this respect the recent works 
\cite{RG1,RG2}. These authors work indeed at the level of the
Saint-Venant equations (i.e. they neglect the non-hydrostatic terms of
the Green-Naghdi equations) but provide a thorough description of
vorticity generation in roll waves and hydraulic jumps for
instance. In their approach, the classical Saint-Venant equations are
coupled with a third equation describing the creation of enstrophy by
wave breaking and by the bottom boundary layer. This enstrophy is
closely related to the tensor $E$ in (\ref{GNvort2dfin}); it seems
therefore possible to combine our approach and Richard and
Gavrilyuk's, leading to a complete Green-Naghdi model describing both
the coupling between surface waves and underlying flows and the
creation of vorticity. This is left for future work.

\vspace{1cm}

Acknowledgment. A.C is support by the grant MTM2011-266696 (Spain), ICMAT Severo Ochoa project SEV-2011-0087 and ERC grant 307179-GFTIPFD. D. L. acknowledges support from the ANR-13-BS01-0003-01
DYFICOLTI, the ANR BOND, and the INSU-CNRS LEFE-MANU project Soli.

\end{document}